\documentclass[12pt]{article}
\usepackage{amssymb}
\input epsf
\def\a{\kern+.6ex\lower.42ex\hbox{$\scriptstyle \iota$}\kern-1.20ex a}
\def\e{\kern+.5ex\lower.42ex\hbox{$\scriptstyle \iota$}\kern-1.10ex e}
\def\n{\'{n}}
\newcommand{\Teis}[2]{
   \setlength{\unitlength}{1ex}
   \begin{picture}(2,0)(0,0.4)
      \put(0,1.1){\line(1,0){2}}
      \put(0,0.9){\line(1,0){2}}
      \put(1,1.2){\makebox(0,0)[b]{$\scriptstyle #1$}}
      \put(1,0.8){\makebox(0,0)[t]{$\scriptstyle #2$}}
   \end{picture}}
\newcommand{\Teisss}[4]{
   \setlength{\unitlength}{1ex}
   \begin{picture}(#3,3)(0,0.4)
      \put(0,1.15){\line(1,0){#3}}
      \put(0,0.85){\line(1,0){#3}}
      \put(#4,1.3){\makebox(0,0)[b]{$#1$}}
      \put(#4,0.7){\makebox(0,0)[t]{$#2$}}
   \end{picture}}
\newcommand{\bbN}{{\mathbb N}}
\newcommand{\bbQ}{{\mathbb Q}}
\newcommand{\bbR}{{\mathbb R}}
\newcommand{\bbC}{{\mathbb C}}
\newcommand{\bbP}{{\mathbb P}}
\newcommand{\bbZ}{{\mathbb Z}}
\newcommand{\bA}{{\mathbf A}}
\newcommand{\bB}{{\mathbf B}}
\newcommand{\bC}{{\mathbf C}}
\newcommand{\bJ}{{\mathbf J}}

\newcommand{\bw}{{\mathbf w}}

\newcommand{\g}{\mathbf g}
\newcommand{\h}{\mathbf h}
\newcommand{\cL}{{\pounds}}
\newcommand{\ccL}{\,\,\tilde{\!\!\pounds}}
\newcommand{\cN}{{\cal N}}

\newcommand{\cT}{{\cal T}}
\newcommand{\oT}{\overline{\cal T}}

\newcommand{\cB}{{\cal B}}

\newcommand{\cZ}{{\cal Z}}
\newcommand{\oE}{{\bar{\cal E}}}
\newcommand{\cE}{{\cal E}}
\newcommand{\oK}{{\overline{\cal K}}}

\newcommand{\nn}{{\bar{n}}}

\def\P{\partial}
\def\GCD{\mbox{\rm GCD}}

\def\cex{\mbox{\rm c.ex.}}
\def\ord1{\mbox{\rm ord}}

\def\incl{\mbox{\rm incl}}
\def\incls{\mbox{\rm\scriptsize incl}}
\def\reds{\mbox{\rm\scriptsize red}}
\def\Track{\mbox{\rm Track}}
\def\Tracks{\mbox{\rm\scriptsize Track}}
\def\fins{\mbox{\rm\scriptsize fin}}
\def\cycle{\mbox{\rm cycle}}
\def\cycles{\mbox{\rm\scriptsize cycle}}
\def\ords{\mbox{\rm\scriptsize ord}}
\def\ints{\mbox{\rm\scriptsize int}}
\def\mins{\mbox{\rm\scriptsize min}}

\def\maxs{\mbox{\rm\scriptsize max}}
\def\grad{\mbox{\rm grad}}

\def\max{\mbox{\rm max}}
\def\sup{\mbox{\rm sup}}
\def\min{\mbox{\rm min}}

\def\supp{\mbox{\rm supp}}
\def\conv{\mbox{\rm conv}}

\def\In{\mbox{\rm in}}            
\def\Zer{\mbox{\rm Zer}}   
\def\char{\mbox{\rm char}}

\def\chars{\mbox{\rm\scriptsize char}}
\newtheorem{thm}{Theorem}[section]
\newtheorem{lm}[thm]{Lemma}
\newtheorem{cor}[thm]{Corollary}
\newtheorem{ex}[thm]{Example}
\newtheorem{prop}[thm]{Proposition}

\newtheorem{rem}[thm]{Remark}
\newtheorem{proper}[thm]{Property}
\newtheorem{df}[thm]{Definition}
\renewcommand{\arraystretch}{0.7}
\title{ON THE \L{}OJASIEWICZ EXPONENT, SPECIAL DIRECTION AND
MAXIMAL POLAR QUOTIENT
\footnotetext{
\begin{minipage}[t]{5.2in}
{\small
{\it Mathematics Subject Classification\/} (2000) 32S55; 14H20;\\
Key words and phrases: plane curve singularity, \L{}ojasiewicz
exponent, polar curve, polar invariants, polar quotients, Eggers tree,
ultrametric space of branches.
Supported in part by the Polish grant: MSHE---No N N201 386634.
} 
\end{minipage}
}}
\begin{document}
\author{Andrzej Lenarcik}
\maketitle
\begin{center}
{\it Professor Stanis\l{}aw \L{}ojasiewicz in memoriam}
\end{center}
\begin{abstract}
\noindent
For a local singular plane curve germ $f(X,Y)=0$
we characterize all nonsingular $\lambda\in\bbC\{X,Y\}$
such that the \L{}ojasiewicz exponent of $\grad\,f$
is not attained on the polar curve $\bJ(\lambda,f)=0$.
When $f$ is not Morse we prove that for the same
$\lambda$'s the maximal polar quotient $q_0(f,\lambda)$
is strictly less than its generic value $q_0(f)$.
Our main tool is the Eggers tree of singularity constructed
as a decorated graph of relations between balls in the space of branches
defined by using a logarithmic distance.
\end{abstract}

\section{Introduction, main results}\label{introduction}

Let $\bbC\{X,Y\}$ be the ring of convergent power series in two variables.
If $f=f_1^{m_1}\dots f_r^{m_r}$ is a decomposition of $f$ into irreducible
pairwise coprime factors in $\bbC\{X,Y\}$ then we put $f_{\reds}=f_1\dots f_r$.
We call $f$ {\it reduced\/} if $f=f_{\reds}$.

For a nonzero series $f=\sum c_{\alpha\beta}X^\alpha Y^\beta\in\bbC\{X,Y\}$
we define the {\it order\/} $\ord1\,f$ as the minimum of $\alpha+\beta$
corresponding to nonzero $c_{\alpha\beta}$ and the {\it initial form\/}
$\In f=\sum_{\alpha+\beta=\ords f}c_{\alpha\beta}X^\alpha Y^\beta$.
We put $\ord1\,0=\infty$ by convention. We call $f$ {\it singular\/}
if $2\leq\ord1\,f<\infty$, {\it nonsingular\/} if $\ord1\,f=1$ and
a {\it unit\/} if $\ord1\,f=0$.

For $f,g\in\bbC\{X,Y\}$ of positive orders we say that
$f$ and $g$ are {\it transverse\/} if the system $\In f=\In g=0$
has no solutions in $\bbC^2\setminus\{0\}$. Otherwise we call
$f$ and $g$ {\it tangent\/}. By $t=t(f)=\ord1(\In\,f)_{\reds}$ we denote
the number of different tangents of $f$.
We call $f$ {\it unitangent\/} if $t(f)=1$ and {\it multitangent\/}
if $t(f)>1$.

Let $f\in\bbC\{X,Y\}$ be a nonzero series without constant term.
The series $f$ defines the curve germ $f=0$ at $0\in\bC^2$. We extend
the term: singular (nonsingular, unitangent, multitangent) for germs and
the term: tranverse (tangent) for pairs of germs.
The singularity $f=0$ is isolated if and only if $f$ is reduced.
Whenever we write a ``singularity'' in this article we mean an ``isolated
singularity''.

Assume that $f$ is reduced ($\ord1\,f\geq 1$).
The {\it\L{}ojasiewicz exponent\/} of $f$ with respect
to a subset $A\subset\bbC^2$, $0\in\overline{A\setminus\{0\}}$, is defined to be
\begin{equation}\label{Lojasiewicz}
  \cL_0(f|A)=\inf\{\theta\geq 0:\;
  |\grad\,f(z)|\geq c|z|^\theta
  \mbox{ for $z\in A$ near zero in $\bbC^2$, $c>0$}\}.
\end{equation}
We write $\cL_0(f)$ for $\cL_0(f|\bbC^2)$.
For nonsingular $f$ we have $\cL_0(f)=0$.
When $\cL_0(f)=\cL_0(f|A)$ we say that the
\L{}ojasiewicz exponent $\cL_0(f)$ {\it is attained\/} on~$A$.
Let $\lambda\in\bbC\{X,Y\}$ be a {\it regular parameter\/} (i.e.
$\lambda(0)=0$, $\lambda$ nonsingular). Consider the germ
$\Gamma_{f,\lambda}$ of {\it polar curve\/}
$$
  \bJ(\lambda,f)=\frac{\P\lambda}{\P X}\frac{\P f}{\P Y}-
               \frac{\P\lambda}{\P Y}\frac{\P f}{\P X}=0\;.
$$
\begin{df}\mbox{}\\{\rm
(a) We define $\lambda$ to be
a {\it special parameter\/} for $f$ if the \L{}ojasiewicz exponent
$\cL_0(f)$ is not attained on $\Gamma_{f,\lambda}$.

\vspace{1ex}\noindent
(b) A~direction $\bw\in\bbP^1(\bbC)$ is defined
to be a {\it special direction\/} of $f$ if there exists a special parameter
$\lambda$ tangent to $\bw$.
} 
\end{df}
One of the goals of this paper is to describe all special parameters
as well as all special directions of $f$. After M.~Lejeune-Jalabert
and B.~Teissier~\cite{LT} we know that for the generic direction
$(a:b)\in\bbP^1(\bbC)$ the parameter $\lambda=bX-aY$ is not special for $f$.
For a~mapping $(f_1,f_2):(\bbC^2,0)\rightarrow(\bbC^2,0)$,
$f_1,f_2\in\bbC\{X,Y\}$, with
isolated zero, the \L{}ojasiewicz exponent can be defined analogously
to~(\ref{Lojasiewicz}). Ch\a{}dzy\n{}ski and Krasi\n{}ski~\cite{ChK} proved
that this exponent is attained on $\{f_1=0\}$ or on $\{f_2=0\}$. This
result applied to the gradient of singularity $f=0$ after coordinate
change can be written as
\begin{thm}\label{thm_ChK1988}{\rm(\cite{ChK}, Main Theorem).}
Let $\lambda,\mu$ be two transversal regular parameters. Then
the \L{}ojasiewicz exponent $\cL_0(f)$ is attained on
$\Gamma_{f,\lambda}$ or on $\Gamma_{f,\mu}$.
\end{thm}
\begin{cor}
A singularity $f=0$ has at most one special direction.
\end{cor}
The following result was obtained independently by
Bogus\l{}awska~\cite{B} and by Kuo and Parusi\n{}ski~\cite{KP1}.
After coordinate change it can be written as
\begin{thm}\label{thm_B_KP}{\rm(\cite{B}, Theorem~2 and \cite{KP1}, Theorem~3.1).}
Let $\lambda$ be a regular parameter transversal to the singularity $f=0$.
Then the \L{}ojasiewicz exponent $\cL_0(f)$ is attained on~$\Gamma_{f,\lambda}$.
\end{thm}
\begin{cor}
If the special direction of $f=0$ exists it is tangent to $f=0$.
\end{cor}
We are going to consider the following problems: (1) to find the conditions
for the existence of the special direction for singularity $f=0$;
(2)~if this direction exists, to determine its position for multitangent $f$;
(3)~to decide: whether or not every regular parameter tangent to the
special direction is special for $f$? Theorem~\ref{main_result}
explains (1) and (2) as well as gives a positive answer to (3).
We call $f=f^{(1)}\dots f^{(t)}$ a {\it tangential decomposition\/} of $f$
if the components $f^{(1)},\dots,f^{(t)}$ are unitangent and pairwise
transverse.
\begin{thm}\label{main_result}{\rm(Main Result A)}\\
Let $f=0$ be a singularity and
let $f=f^{(1)}\dots f^{(t)}$ be a tangential decomposition of $f$ $(t\geq 1)$.
Then
\begin{itemize}
\item[\rm(i)]
$\displaystyle\cL_0(f)=%
\begin{array}[t]{c}\max\\{\scriptstyle i=1,\dots,t}\end{array}%
(\cL_0(f^{(i)})+\ord1\,f-\ord1\,f^{(i)})\;.$
\item[\rm(ii)] Let $\lambda$ be a regular parameter.
If the maximum in {\rm(i)} is realized for exactly one
index $i_0\in\{1,\dots,t\}$ then $\lambda$ is special
for $f$ if and only if $\lambda$ is tangent to $f^{(i_0)}$.
\item[\rm(iii)] If the maximum in {\rm(i)} is realized for two or more
indicies from $\{1,\dots,t\}$ then there are no special parameters
for $f$.
\end{itemize}
\end{thm}
We prove (i) of Theorem~\ref{main_result} in
Section~\ref{Eggers_tree} and (ii), (iii) in Section~\ref{factorization}.

\vspace{1ex}\noindent
When all tangential components
of $f$ are nonsingular, we call $f=0$ an {\it ordinary singularity\/}.
If additionally $\ord1\,f=2$ then we call $f=0$ a {\it Morse singularity\/}.
\begin{cor}\label{cor12}
Assume that $f=0$ is an ordinary singularity. Then for every local
parameter $\lambda$
$$
  \cL_0(f)=\cL_0(f|\Gamma_{f,\lambda})=\ord1\,f-1\;.
$$
\end{cor}
\begin{cor}\label{cor}
The tangent direction of any unitangent singularity is special.
\end{cor}
\begin{ex}\label{ex14}{\rm
Let $f=f^{(1)}f^{(2)}$ where $f^{(1)}=Y^5+X^2$ and $f^{(2)}=Y(Y^2-X^4)$.
By direct computation (or for example by using~\cite{L1}) we obtain
$\cL_0(f^{(1)})=4$ and $\cL_0(f^{(2)})=5$. We have
$\cL_0(f^{(i)})+\ord1\,f-\ord1\,f^{(i)}=7$ for $i=1,2$. By
Theorem~\ref{main_result} $\cL_0(f)=7$ and the special direction does not
exist.
} %
\end{ex}
\begin{rem}\label{rem_gwoz}{\rm
An interesting family of examples of singularities without special
direction was proposed by Gwo\'{z}dziewicz (oral communication).
Let $f\in\bbC\{X,Y\}$ be such that $f(X,Y)=f(Y,X)$ with the only
tangents $X=0$ and $Y=0$. By symmetry of $f$ and both
Theorems~\ref{thm_ChK1988} and~\ref{thm_B_KP}
we conclude that there are no special parameters.
We do not need Theorem~\ref{main_result}.
} %
\end{rem}
Let us recall some facts concerning the \L{}ojasiewicz exponent $\cL_0(f)$
of a holomorphic function defined by $f\in\bbC\{X_1,\dots,X_n\}$
with an isolated singularity at zero. Let $[x]$ stand for the integer
part of $x$. Lu and~Chang~\cite{LuChang} (developing the results
of Kuo~\cite{Kuo1969}, Kuiper~\cite{Kuiper}, Bochnak and
\L{}ojasiewicz~\cite{BochLoj}) proved that adding to $f$ monomials of order
greater than $[\cL_0(f)]+1$ does not change the topological type of
singularity $f=0$.
The minimal integer with this property is called the $C^0$-sufficiency degree
of $f$. Teissier~\cite{T1} showed that this degree equals
$[\cL_0(f)]+1$ (Kucharz~\cite{Kuch} found an example
that the analogous equality is not true in the real case).
In the same paper Teissier found a relation between the \L{}ojasiewicz
exponent and the {\it maximal polar invariant\/}.
References to papers concerning the different kinds of the
\L{}ojasiewicz exponents can be found in~\cite{RodSpodz}.

In dimention two Kuo and Lu~\cite{KuoLu} described
$\cL_0(f)$ in terms of a tree model constructed on the basis of Puiseux roots
of $f=0$. Following Teissier's result, authors focused their attention on polar
invarians and so called {\it polar quotients\/}. A survey of results
concerning this subject in dimension two is given in~\cite{GLP1}.
We explain the notions of polar quotients and polar invariants for curve germs.
For any $f,g\in\bbC\{X,Y\}$ the {\it intersection multiplicity\/}
$(f,g)_0$ is defined to be the $\bbC$-codimension of the ideal generated
by $f$ and $g$ in $\bbC\{X,Y\}$.
Take an irreducible $h\in\bbC\{X,Y\}$. We call $h$, as well as
the corresponding germ $h=0$, a {\it branch\/} ({\it smooth\/} branch
if $h$ is nonsingular). The {\it semigroup\/} of $h$ is
\begin{equation}\label{semigroup}
\Gamma(h)=\{(h,g)_0:\,g\in\bbC\{X,Y\}, h\mbox{ does not divide }g\}\;.
\end{equation}
Now, let $f,g\in\bbC\{X,Y\}$ be reduced series. We call two germs
$f=0$ and $g=0$ {\it equisingular\/} if there exist factorizations
$f=f_1\dots f_r$ and $g=g_1\dots g_s$ into branches such that
$r=s$, $\Gamma(f_i)=\Gamma(g_i)$ for $i=1,\dots,r$
and $(f_i,f_j)_0=(g_i,g_j)_0$ for $i,j=1,\dots,r$. Equisingularity relation
defines {\it equisingularity classes\/} in the set of germs.
By an {\it equisingularity invariant\/} we mean a function constant
in every equisingularity class.

For $f,\lambda\in\bbC\{X,Y\}$
($f$ singular, $\lambda$ nonsingular) let us consider the set of
{\it polar quotients\/} of $f$ with respect to parameter $\lambda$:
\begin{equation}\label{quotients}
Q(f,\lambda)=\left\{\frac{(f,g)_0}{(\lambda,g)_0}:\;
g\mbox{ irreducible factor of }
\bJ(\lambda,f),\,g\neq\lambda\right\}\;.
\end{equation}
We define the {\it maximal polar quotient\/}
$q_0(f,\lambda)$ as $\max\,Q(f,\lambda)$ if $Q(f,\lambda)\neq\varnothing$ and
as $-\infty$, otherwise (for the case $Q(f,\lambda)=\varnothing$ see Example~%
\ref{xy} and Remark~\ref{Q_empty}). Teissier proved that the set 
$Q(f):=Q(f,\lambda)$ does not depend on sufficiently generic $\lambda$ 
and that it is an equisingularity invariant of $f$. 
We call $Q(f)$ the set of {\it polar invariants\/}. It is always
nonempty for singular $f$. Then $q_0(f)=\max\,Q(f)$ is called the
{\it maximal polar invariant\/}.  Teissier \cite{T1} proved that
\begin{equation}\label{Te}
\cL_0(f)=q_0(f)-1\;.
\end{equation}
Analogously, as we did for the germs, we define the {\it equisingularity of pairs\/}
$(f,\lambda)$, $(g,\mu)$~\cite{GP3}. We consider
{\it equisingularity classes\/} and {\it equisingularity invariants\/} for
pairs. According to~\cite{GP2002} we know that the set $Q(f,\lambda)$ is
an invariant in this sense. Assume that $f,\lambda$ are transverse. In this
case P\l{}oski \cite{P4} showed $q_0(f)=q_0(f,\lambda)$ (the equality
$Q(f)=Q(f,\lambda)$ was shown in \cite{GP2002}); by Theorem~\ref{thm_B_KP}
and (\ref{Te}) we obtain $\cL_0(f)=\cL_0(f|\Gamma_{f,\lambda})=q_0(f,\lambda)-1$.
In the following theorem we explain relations
between these numbers for an arbitrary $\lambda$.
\begin{thm}\label{thm_Loj_and_polar}{\rm(Main Result B)}\\
Let $f=0$ be a singularity and let $\lambda\in\bbC\{X,Y\}$ be a regular
parameter. Then:
\begin{itemize}
\item[\rm(a)]
$\cL_0(f)\geq\cL_0(f|\Gamma_{f,\lambda})\geq q_0(f,\lambda)-1$.
\item[\rm(b)] Moreover, if $f=0$ is not Morse then the equalities
$\cL_0(f)=\cL_0(f|\Gamma_{f,\lambda})$ and $q_0(f)=q_0(f,\lambda)$ are
satisfied for exactly the same $\lambda$'s.
\end{itemize}
\end{thm}
We prove this theorem in Section~\ref{factorization}.
\begin{ex}\label{xy}{\rm
Assume that $f=XY$ and $\lambda=X$. By direct computation we obtain
$\cL_0(f)=\cL_0(f|\Gamma_{f,\lambda})=1$ and $q_0(f)=2$. But
$\bJ(\lambda,f)=(\P f/\P Y)=X$. Hence $Q(f,\lambda)=\varnothing$ and
$q_0(f)>q_0(f,\lambda)=-\infty$. This explains the assumption
``$f$ is not Morse'' in Theorem~\ref{thm_Loj_and_polar}(b)
and later in Corollary~\ref{cor18}.
} %
\end{ex}
Let us observe two corollaries of Theorem~\ref{thm_Loj_and_polar}.
The first one is strightforward.
\begin{cor}\label{cor17} {\rm(see~\cite{P4}, Corollary~1.4)}\\
For a singularity $f=0$ and a regular parameter $\lambda$ we have
$q_0(f)\geq q_0(f,\lambda)$.
\end{cor}
As a consequence we obtain
\begin{cor}\label{cor18}
Let $\lambda,\mu$ be two transversal regular parameters.
Then if $f$ is not Morse
then
$$
q_0(f)=\max\{q_0(f,\lambda),q_0(f,\mu)\}\;.
$$
\end{cor}
Proof. From the quoted result of Ch\a{}dzy\n{}ski and Krasi\n{}ski we can
assume that $\cL_0(f)=\cL_0(f|\Gamma_{f,\lambda})$. Hence by Theorem~\ref{thm_Loj_and_polar}(b)
$q_0(f)=q_0(f,\lambda)$.
We finish the proof by using Corollary~\ref{cor17}~\rule{1ex}{1ex}

Our main tool is the Eggers tree~\cite{E,GB,Wall1,Wall2,McNeal-Nemethi}
which is a decorated graph that represents the equisingularity class of
a germ $f=0$.

In Section~\ref{Eggers_tree} we propose a new construction of the Eggers
tree of $f=0$ by using the {\it order of contact\/} of P\l{}oski~\cite{P1}.
We do not need Puisex series which were used in the original construction~\cite{E,GB}.
P\l{}oski proved that the order of contact of every two branches
satisfies the axioms of logarithmic distance. This distance allows us
to define {\it characteristic contacts\/} (\ref{contact}) for every
singular branch. We can also consider balls (every branch inside the ball is
a center of this ball). We assign to the germ $f=0$
the set of balls called {\it Eggers collection\/} (Definition~\ref{collection}).
In this collection we have the balls that come from intersections of branches
and the balls that come from singular branches and their characteristic
contacts. The {\it Eggers tree\/} is a graph determined by the Eggers
collection (Definition~\ref{df_Eggers_tree}).
The balls correspond to vertices of the graph. The edges correspond to
inclusions of successive balls.

It is recently proved~\cite{GP2011} that the order of contact
satisfies the axioms of logarithmic distance also in
positive characteristic. This suggests an application of this new construction
for singularities over an arbitrary field.

As an application of the Eggers tree technique we give a recursive version
of Eggers formula for polar invariants $Q(f)$ (\ref{q}).
The formula for $Q(f)$ together with (\ref{Te}) suffices to prove
Theorem~\ref{main_result}~(i).

In order to describe a position of an arbitrary branch $h$ with
respect to the germ $f=0$ we consider the ball $B_f(h)$ with $h$ as a center.
The radius of $B_f(h)$ equals the maximal order of contact of $h$ with branches
of $f=0$ (see: a definition before Property~\ref{proper_extension}).

Let $\lambda$ be an arbitrary regular parameter (possibly a branch of the
germ $f=0$). In Section~\ref{Loj_rel} we give formulas for the \L{}ojasiewicz
exponent $\cL_0(f|\Gamma_{f,\lambda})$ (Proposition~\ref{prop_Loj_rel},
Corollary~\ref{cor_Loj_rel}). These formulas involve the position of
$\lambda$ as well as the positions of branches of the polar $\bJ(\lambda,f)$
with respect to $f=0$. We show (Example~\ref{ex_Loj_irr})
equisingular pairs $(f,\lambda)$, $(f',\lambda')$ such that
$\cL_0(f|\Gamma_{f,\lambda})\neq\cL_0(f'|\Gamma_{f',\lambda'})$.
Hence $\cL_0(f|\Gamma_{f,\lambda})$ is not in general an equisingularity
invariant of the pair $(f,\lambda)$. This example concerns the very specific
equisingularity class when $f=0$ is unitangent and $B_f(\lambda)$ coincides
with the unique ball of the Eggers collection. For each different class
the \L{}ojasiewicz exponent
$\cL_0(f|\Gamma_{f,\lambda})$ is an invariant (see: Lemma~\ref{lm_Loj_rel}).

In Section~\ref{factorization} we propose Theorem~\ref{thm_factorization}
to factorize the polar $\bJ(\lambda,f)$ involving only the equisingularity
information of the pair $(f,\lambda)$. If $\lambda$ is tranversal to $f=0$
then for every factor $g$ of $\bJ(\lambda,f)$ the ball $B_f(g)$ belongs to
the Eggers collection. When $\lambda$ is tangent to $f=0$ the position of
$B_f(g)$ in the Eggers collection is not in general determined by the
eqisingularity class of the pair $(f,\lambda)$. This phenomenon was observed
by Kuo and Parusi\'{n}ski~(\cite{KP2}, Example~8.1) for Kuo-Lu trees~\cite{KuoLu}.
In this case we assign $g$ to the nearest succesive ball in the Eggers collection.
Finally, every factor $g$ of $\bJ(\lambda,f)$ (different from $\lambda$)
is assigned to a ball $B$ of the Eggers tree or to the ball $B=B_f(\lambda)$.
The ``packages'' of $g$'s form factors $h_B$ of $\bJ(\lambda,f)$.
In Theorem~\ref{thm_factorization}~(ii) we desribe the contacts of $g's$
with $\lambda$. In Theorem~\ref{thm_factorization}~(iii) we give two formulas
for $(h_B,\lambda)_0$. The first one is analogous to that from~\cite{E,GB};
the second concerns the ball $B_f(\lambda)$ and the balls from the Eggers
collection which have $\lambda$ as their centre. As a consequence of
Theorem~\ref{thm_factorization} we obtain a version of the result of Eggers
(Corollary~\ref{cor_Eggers}). For typical equisingularity classes,
different from the class of Example~\ref{ex_Loj_irr},
we describe $\cL_0(f|\Gamma_{f,\lambda})$ in Lemma~\ref{lm_Loj_rel}.
This lemma allows us to prove Theorem~\ref{main_result}~(ii,iii).
We obtain formulas for polar quotients~(Proposition~\ref{prop_quotients}),
for their multiplicities (Remark~\ref{rem_mult}) and
for the maximal polar quotient $q_0(f,\lambda)$ (Lemma~\ref{polar_and_tree}).
Applying Lemmas~\ref{lm_Loj_rel} and~\ref{polar_and_tree}
we prove Theorem~\ref{thm_Loj_and_polar}.

Theorem~\ref{thm_factorization} is a version of known results
(\cite{E,LMW1,LMW2,GB,GP2002,Ma,KP2}).
We generalize~\cite{E,LMW1,LMW2,GB,GP2002}.
In~\cite{E,GB} (resp. in~\cite{LMW1,LMW2}) $\lambda$ is generic
(resp. $\lambda$ is transversal to $f$) whereas in
Theorem~\ref{thm_factorization} $\lambda$ is an arbitrary regular
parameter. In comparison to~\cite{GP2002}, where $Q(f,\lambda)$ is described
in terms of the equisingularity class of the pair $(f,\lambda)$, we give
formulas for the multiplicities of polar quotients.
The paper of Maugendre~\cite{Ma} concerns a more general situation of
{\it jacobian quotients\/}. For nonzero series $f,g\in\bC\{X,Y\}$ without
constant terms the {\it jacobian curve\/} $\bJ(f,g)=0$ is considered.
Every branch $h$ of $\bJ(f,g)$ which is not a branch of $fg$ defines
a jacobian quotient $(f,h)_0/(g,h)_0$. Maugendre described the set of
jacobian quotients in terms of the minimal resolution of $f\cdot g$.
Applying this result with smooth $g$ we can obtain the set of polar quotients
but without multiplicities. Kuo and Parusi\n{}ski~\cite{KP2} considered
the case when the Puiseux roots of $fg$ are different. They constructed
a tree model $T(f,g)$ similar to that of~\cite{KuoLu}. They described
how the Puiseux roots of $\bJ(f,g)$ ``leave'' $T(f,g)$. This construction
depends on the choice of the coordinate system. 
It is possible to apply this result
to prove Theorem~\ref{thm_factorization}, but it requires effort
to move from Puiseux roots to branches and to eliminate an influence of
the coordinate system. Finally, we decided to present in
Section~\ref{factorization_proof} a self-contained proof based
on the technique of paths of the Newton algorithm from~\cite{Len2004}.

From Theorem~\ref{thm_factorization} it follows that the polar quotients
together with their multiplicities are equisingularity invariants of the
pair: germ, regular parameter (see~\cite{GLP2}). The analogous fact for
jacobian pairs was recently proved in~\cite{Mi,Gwoz}.

\section{The Eggers tree}\label{Eggers_tree}

In this section we construct the Eggers tree by using the order
of contact of P\l{}oski. We propose a recursive version of the Eggers
formula for the polar invariants (\ref{Qf}). By using the formula
we prove Theorem~\ref{main_result}(i).

Let us denote by $\cB$ the set of all branches.
From P\l{}oski~\cite{P1,ChP} we know that for branches $f,g\in\cB$
the {\it order of contact\/}
$$
  d(f,g)=\frac{(f,g)_0}{(\ord1\,f)(\ord1\,g)}
$$
satisfies the axioms of {\it logarithmic distance\/}:
\begin{itemize}
\item[$(D_1)$] $d(f,g)=\infty$ if and only if the germs $f=0$ and $g=0$ coincide,
\item[$(D_2)$] $d(f,g)=d(g,f)$,
\item[$(D_3)$] $d(f,g)\geq\min\{d(f,h),\,d(g,h)\}$.
\end{itemize}
Since $(f,g)_0\geq(\ord1\,f)(\ord1\,g)$ we have $d(f,g)\geq 1$.
Moreover $d(f,g)=1$ if and only if $f$, $g$ are transverse.
A simple consequence of ($D_3$) is
\begin{itemize}
\item[$(D_3')$] If $d(f,h)\neq d(g,h)$ then $d(f,g)=\min\{d(f,h),\,d(g,h)\}$.
\end{itemize}

\vspace{1ex}\noindent{\bf Characteristic contacts}\\
Recall that the semigroup of a branch $f$ can be written as
$\Gamma(f)=\bbN\bar{\beta}_0+\dots+\bbN\bar{\beta}_{\g}$ where
$\bar{\beta}_0<\dots<\bar{\beta}_{\g}$ is the minimal sequence of
{\it semigroup generators\/}. We call $\g=\g(f)$ the
number of {\it characteristic pairs\/} of $f$.
For smooth branches we have $\g=0$, $\bar{\beta}_0=1$.
For $k=1,2,\dots$ we define the {\it characteristic contacts\/}
\cite{LJ,GB-P}
\begin{equation}\label{contact}
   d_k(f)=\sup\{d(f,h):\,h\in\cB,\,\g(h)<k\}\;.
\end{equation}
For $k>\g(f)$ we have $d_k(f)=\infty$. For singular branch $f$ we have
\begin{equation}
   d_k(f)=\frac{\GCD(\bar{\beta}_0,\dots,\bar{\beta}_{k-1})\bar{\beta}_k}
   {(\bar{\beta_0})^2}\quad\mbox{for }k=1,\dots,\g(f)\;.
\end{equation}
We have $d_1<\dots<d_\g$ which is equivalent to
$n_k\bar{\beta}_k<\bar{\beta}_{k+1}$ ($k=1,...,\g(f)-1$).
We write $\char(f)=\{d_1,\dots,d_{\g}\}$.
By $(n_1,\dots,n_g)$ we denote the corresponding sequence
$n_k:=\GCD(\bar{\beta}_0,\dots,\bar{\beta}_{k-1})/%
\GCD(\bar{\beta}_0,\dots,\bar{\beta}_k)$, $k=1,\dots,\g$. We have
(compare~\cite{Ma}, Proposition~3.2)
\begin{equation}\label{charcontacts}
d_1\in\frac{\bbN}{n_1}\setminus \bbN\mbox{ and }
d_k\in\frac{\bbN}{(n_1\dots n_{k-1})^2n_k}\setminus
\frac{\bbN}{(n_1\dots n_{k-1})^2},\quad k=2,\dots,\g\;.
\end{equation}
Let us denote $\nu_0=1,\nu_1=n_1,\dots,\nu_{\g}=n_1\dots n_{\g}$
($\nu_{\g}=\ord1\,f$). The formula
\begin{equation}\label{n}
  n_k=\min\left\{n\geq 1:\,d_k\in\frac{\bbN}{\nu_{k-1}^2n}\right\},\quad
  k=1,\dots,\g
\end{equation}
enables us to reconstruct the sequence $(n_1,\dots,n_{\g})$ from
$(d_1,\dots,d_{\g})$.

The following classical facts are useful.
\begin{proper}\label{proper_useful}\mbox{}
\begin{itemize}
\item[\rm 1)] For $f\in\cB$, $R\in\bbQ\cap\langle 1,\infty)$ there exists
$g\in\cB$ such that $d(f,g)=R$.
\item[\rm 2)]
For $f\in\cB$, $\g=\g(f)>0$ there exists a sequence of branches
$f_0,\dots,f_{\g-1}$ such that $\g(f_k)=k$ and $d(f,f_k)=d_{k+1}(f)$
for $k=0,\dots,\g-1$.
\end{itemize}
\end{proper}
For singular $f$ in 2) $f_0$ is the classical maximal contact of Hironaka.

\vspace{1ex}\noindent{\bf Balls and trees}\\
Let $f\in\cB$ and let $R\in\langle 1,\infty\rangle$.
The set $B(f,R)=\{g\in\cB:\,d(f,g)\geq R\}$
will be referred to as the {\it ball\/} with center $f$ and radius $R$.
By using ($D_3$) we can prove that every element of the ball is a center of this ball.
Clearly $\cB=B(f,1)$ for $f\in\cB$. For $f,g\in\cB$ we put $B(f,g)=B(f,d(f,g))$.
For each ball $B$ we define the {\it diameter\/} $d(B)=\inf\{d(f,g):\,f,g\in B\}$
which is equal to the radius. For any two balls $B,B'$ if
$B\cap B'\neq\varnothing$ then $B\subset B'$ or $B\supset B'$.
We define $B\leq B'$ if $B\supset B'$ and $B<B'$ if $B\leq B'$ and $B\neq B'$.
Let $f\in\cB$ and let $R,R'\geq 1$. By Property~\ref{proper_useful}-1 we obtain
\begin{equation}
B(f,R)=B(f,R')\Leftrightarrow R=R'\;.
\end{equation}
Now we want to define the {\it Eggers\/} collection of a singularity. It is
a finite set of balls. Let us consider a germ $f=0$ and the
factorization $f=f_1\dots f_r$ into branches;
$r=r(f)$ is the number of branches ($r\geq 1$).
\begin{df}\label{collection} {\rm(Eggers collection)\\
By the {\it Eggers collection\/} of the germ $f=0$ we mean
the collection of balls
$$
\oE(f)=\{B(f_i,f_j)\}_{i,j=1,\dots,r}\cup
\bigcup_{\begin{array}{c}\scriptstyle i=1,\dots,r\\
\scriptstyle f_i\,\mbox{\scriptsize singular}\end{array}}
\{B(f_i,d_{i,k})\}_{k=1,\dots,g(f_i)}
$$
where $\{d_{i,1},\dots,d_{i,\g(f_i)}\}$ are
the characteristic contacts of singular branches.
} %
\end{df}
Let us observe that the balls $B(f_i,f_i)$, $i=1,\dots,r$, of infinite
diameters are in the collection. Balls of finite diameters form the
{\it truncated Eggers collection\/} $\cE(f)$.
For a smooth branch $f$ we have $\oE(f)=\{B(f,f)\}$ and $\cE(f)=\varnothing$.
The following proposition is a consequence of Property~\ref{proper_useful}.
\begin{prop}
Let us consider a ball $B$.
\begin{itemize}
\item[\rm(a)] Then the characteristic contacts
strictly less than $d(B)$ of every branch of $B$ are the same.
\item[\rm(b)] For every $d\geq d(B)$ there exists $f\in B$ such that
$d\notin\char(f)$.
\end{itemize}
\end{prop}
This allows us to define the {\it characteristic of a ball\/} $B$ as
\begin{equation}\label{char_ball}
\char(B)=\bigcap_{f\in B}\char(f)\;.
\end{equation}
\begin{cor}\label{cor_less}
All the elements of $\char(B)$ are strictly less than $d(B)$.
\end{cor}
\begin{cor}\label{cor_smooth_centre}
If $B$ has a smooth centre then $\char(B)=\varnothing$.
\end{cor}
Now, we want to define the number $\nu(B)$ for every ball $B$.
If $\char(B)=\varnothing$ then we put $\nu(B)=1$.
If $\char(B)=(d_1,\dots,d_k)$ with the corresponding
sequence $(n_1,\dots,n_k)$ then we put $\nu(B)=n_1\dots n_k$.
In analogy to (\ref{n}) we define 
\begin{equation}\label{nb}
n(B)=\min\left\{n\geq 1:\,d(B)\in\frac{\bbN}{\nu(B)^2n}\right\}\;.
\end{equation}
We call $B$ a {\it characteristic ball\/} if $n(B)>1$ and
a {\it noncharacteristic ball\/} if $n(B)=1$.
Let B be a ball and let $\cZ$ be a set of balls. We call $B'\in\cZ$
a {\it direct successor\/} of $B$ in $\cZ$ if $B<B'$ and from $B<B_1\leq B'$,
$B_1\in\cZ$, it follows that $B_1=B'$. If additionally $B\in\cZ$ then we call
$(B,B')$ a {\it pair of successive balls\/} in $\cZ$.
\begin{df}\label{df_Eggers_tree}{\rm (Eggers tree)\\
We define the {\it Eggers tree\/} of the germ $f=0$ as a graph whose vertices are
the balls of the Eggers collection $\oE(f)$ and the edges are the pairs
of successive balls $(B,B')$ in $\oE(f)$. From the axioms
($D_1$--$D_3$) it follows that this graph is a rooted tree where the
root is the ball with the minimal diameter. Black (white) vertices
are the balls of finite (infinite) diameters. We call
an edge $(B,B')$ {\it discontinuous\/} when $d(B)\notin\char(B')$ and
{\it solid\/} when $d(B)\in\char(B')$. Moreover, we assign $d(B)$ as a
decoration to every black vertex. By the {\it truncated Eggers tree\/}
we mean the analogous graph constructed for the truncated Eggers
collection $\cE(f)$.
} %
\end{df}
In what follows we denote an edge $(B,B')$ as $B<B'$. By $B_{\mins}(f)$
we denote the ball with the minimal diameter in $\oE(f)$.
\begin{rem}{\rm
Eggers originally assigned {\it contact exponents\/} $\cex(B)$ to every black
vertex of the tree. We can obtain these exponents by the following classical computation.
If $\char(B)=\varnothing$ then $\cex(B)=d(B)$.
If $\char(B)=(d_1,\dots,d_k)$ with the corresponding sequence $(n_1,\dots,n_k)$
then we have
\begin{equation}\label{eq_cex}
\cex(B)=d_1+n_1(d_2-d_1)+\dots+n_1\dots n_k(d(B)-d_k)\;.
\end{equation}
} 
\end{rem}
The following simple property is useful for constructing the tree.
\begin{proper}\label{proper_useful_balls}
Let $B_1=B(f_1,R_1)$ and $B_2=B(f_2,R_2)$. Then
$$
B_1\cap B_2=\varnothing\Leftrightarrow\min(R_1,R_2)>d(f_1,f_2)\;.
$$
\end{proper}
\begin{cor}\label{cor_useful_balls}
For $B_1,B_2$ as above:
\begin{itemize}
\item[\rm(a)] $B_1<B_2\Leftrightarrow\left(\rule{0ex}{2.2ex}R_1<R_2%
\mbox{ and }R_1\leq d(f_1,f_2)\right)$,
\item[\rm(b)] $B_1=B_2\Leftrightarrow R_1=R_2\leq d(f_1,f_2)$.
\end{itemize}
\end{cor}
\begin{ex}\label{ex26}{\rm
Let us consider $f=(Y^5+X^2)Y(Y^2-X^4)$ as in Example~\ref{ex14}
with four irreducible branches
$f_1=Y^5+X^2$, $f_2=Y$, $f_3=Y-X^2$, $f_4=Y+X^2$. The contacts $d(f_i,f_j)$
($i,j=1,\dots,4$) of different branches are presented in the table
$$
\epsfbox{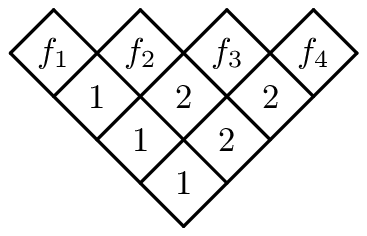}
$$
By using Property~\ref{proper_useful_balls} and Corolary~\ref{cor_useful_balls}
we recognize that the Eggers collection
$\oE(f)$ has three balls with finite diameters
$\cB=B(f_1,f_2)=B(f_1,f_3)=B(f_1,f_4)$,
$B_1=B(f_1,5/2)$, $B_2=B(f_2,f_3)=B(f_2,f_4)=B(f_3,f_4)$
and four balls with infinite diameters
$\{f_1\}$, $\{f_2\}$, $\{f_3\}$, $\{f_4\}$
that can be identified with the branches.
Let us notice that only $f_1$ is singular with $\char(f_1)=\{5/2\}$. Other
branches are smooth. There is only one solid edge $B_1<\{f_1\}$.
All the other edges are discontinuous. The ball $B_1$ is characteristic
whereas $B_2$ and $\cB$ are noncharacteristic.
} %
$$
\epsfbox{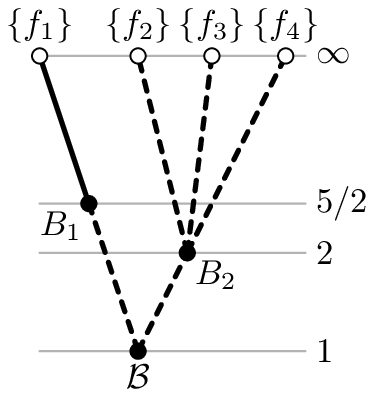}
$$
\end{ex}
Now, let us consider an arbitrary ball $B$ and the set of branches
$\cB_f=\{f_1,\dots,f_r\}$ of the germ $f=0$. By $t_f(B)$ we denote the
number of direct successors $\{B_1,\dots,B_t\}$ of $B$ in $\oE(f)$. By
$t^{(1)}_f(B)$ (resp. $t^{(2)}_f(B)$) we denote the number of direct
successors $B_l$, $l\in\{1,\dots,t_f(B)\}$, that $d(B)\notin\char(B_l)$ (resp.
$d(B)\in\char(B_l)$). Clearly, $t_f(B)=t^{(1)}_f(B)+t^{(2)}_f(B)$. If
$B\in\cE(f)$ then $t^{(1)}_f(B)$ (resp. $t^{(2)}_f(B)$) equals the number of
discontinuous (resp. solid) edges that leave $B$. We have
\begin{prop}\label{prop_t}\mbox{}\\
Let $f_i\sim f_j\Leftrightarrow d(f_i,f_j)>d(B)$ be a relation in the set
$\cB_f\cap B$. Then
\begin{itemize}
\item[\rm(a)] it is an equivalency relation,
\item[\rm(b)] $t_f(B)$ equals the number of equivalency classes of
the relation in $\cB_f\cap B$,
\item[\rm(c)] if $B$ is characteristic then $t^{(1)}_f(B)\in\{0,1\}$.
\end{itemize}
\end{prop}
Proof. (a) is a direct consequence of the axioms ($D_1$-$D_3$). To prove (b)
we first observe that the number of direct successors of $B$ does not
change when we substitute $\oE(f)$ by
$\oE_{\ints}(f):=\{B(f_i,f_j)\}_{i,j=1,\dots,r}$. Then we use the axioms.
To prove (c) let us consider a characteristic ball $B$. It suffices to show
that if $f,g\in B$ and $d(B)\notin\char(f)$, $d(B)\notin\char(g)$ then
$d(f,g)>d(B)$. By Property~\ref{proper_useful}-2 we choose $f',g'\in\cB$
such that $\char(f')=\char(g')$ and $d(f,f')>d(B)$, $d(g,g')>d(B)$.
Clearly $\ord1\,f'=\ord1\,g'=\nu(B)$. Since $d(f,g)\geq d(B)$ then
$d(f',g')\geq d(B)$. We have $d(f',g')\in\bbN/\nu(B)^2$ wherease
$d(B)\in\bbN/(\nu(B)^2n(B))\setminus\bbN/\nu(B)^2$ by~(\ref{charcontacts}).
Hence $d(f',g')>d(B)$ and $d(f,g)\geq\min\{d(f,f'),d(f',g'),d(g',g)\}>d(B)$~%
\rule{1ex}{1ex}

\vspace{1ex}\noindent Let us observe that if $B<B'$ then by (\ref{char_ball})
$\char(B)\subset\char(B')$.
\begin{proper}\label{proper_precise}
Let $B<B'$ be an edge.
Then
\begin{itemize}
\item[\rm(a)]
the edge is discontinuous if and only if $\char(B')=\char(B)$.
\item[\rm(b)] For solid edge $\char(B')=\char(B)\cup\{d(B)\}$.
\end{itemize}
\end{proper}
By a {\it chain\/}
in the Eggers collection (tree) we mean an
increasing sequence of successive balls (vertices).
\begin{rem}{\rm
The equisingularity class of a singularity can be reconstructed
from its Eggers tree.
The branches correspond to white vertices. In order to recognize
the characteristic of a branch we consider the chain that joins
the minimal vertex with the corresponding white vertex and 
we apply Property~\ref{proper_precise}. The contact
between branches $f_i,f_j$ is 
$$
  d(f_i,f_j)=\max\{d(B):\,B\in\oE(f),B\leq\{f_i\},B\leq\{f_j\}\}\;.
$$
} 
\end{rem}

\vspace{1ex}\noindent{\bf Tangential decomposition}\\
Consider the germ $f=0$ with branches $\cB_f=\{f_1,\dots,f_r\}$.
Applying Proposition~\ref{prop_t}(a,b) with $B=\cB$
we divide $\cB_f$ due to the equivalency relation $d(f_i,f_j)>1$.
When we multiply the branches inside each class we obtain a tangential
decomposition $f=f^{(1)}\dots f^{(t)}$ (as in Introduction) where $t=t(f)$
is the number of tangents of the germ. The following property follows directly
from Definition~\ref{collection}.
\begin{proper}\label{proper_tangential}
If $t=t(f)>1$ then $\oE(f)=\cB\cup\oE(f^{(1)})\cup\dots\cup\oE(f^{(t)})$.
\end{proper}

\vspace{1ex}\noindent
{\bf Orders, polar invariants and multiplicities of balls}\\
For an arbirtary ball $B$ and the germ $f=0$ with branches
$\cB_f=\{f_1,\dots,f_r\}$ we define the {\it order\/}
$O_f(B)=\sum_i\ord1\,f_i$ where the summation runs over $f_i\in B$.
It is convenient to define the {\it family of balls determined by\/} $f$:
\begin{equation}\label{def_T}
\oT(f)=\{B\mbox{ ball}:\,O_f(B)>0\}\;.
\end{equation}
We have $\oE(f)\subset\oT(f)$. We write $\cT(f)$ when we omit the balls with
infinite diameters. We say that a ball $B\in\oT(f)$ {\it lies on the edge\/}
$B_1<B_2$ ($B_1,B_2\in\oE(f)$) if $B_1<B\leq B_2$. We define the pair
$(\cB,B_{\mins}(f))$ to be the {\it trunk\/} of $f$. We say that $B$ lies
on the trunk if $\cB\leq B\leq B_{\mins}(f)$. The family $\oT(f)$ contains
exactly these balls that lie on the edges or on the trunk. Let us observe
that the function $B\mapsto O_f(B)$ is constant for balls $B$ lying on the
one edge or on the trunk. This function has ``jumps'' only for
$B\in\cE_{\ints}(f)$. Let us observe that $O_f(B)=\ord1\,f$ for
balls lying on the trunk.

For every ball $B\in\cT(f)$ we define the number $q_f(B)$.
First we define $q_f(B)$ for $B\in\cE(f)$. We consider the unique chain
$B_{\mins}(f)=B_1<B_2<\dots<B_l=B$ and we put
\begin{equation}\label{q}
\renewcommand{\arraystretch}{1.2}
q_f(B_k)=\left\{\begin{array}{ll}
O_f(B_1)d(B_1),&\quad k=1\\
q_f(B_{k-1})+O_f(B_k)(d(B_k)-d(B_{k-1})),&\quad 1<k\leq l\;.\end{array}\right.
\end{equation}
Then we use a linear interpolation due to $d(B)$ to define $q_f(B)$ for
every $B\in\cT(f)$ (i.e if $B$ lies on the edge $B_1<B_2$ then
$q_f(B)=q_f(B_1)+O_f(B)(d(B)-d(B_1))$, clearly $O_f(B)=O_f(B_2)$).
We always have $q_f(\cB)=\ord1\,f$. Let us observe that
\begin{equation}\label{monotonicity1}
\mbox{for }B,B'\in\cT(f)\mbox{ if }B<B'\mbox{ then }q_f(B)<q_f(B')\;.
\end{equation}
To every ball $B\in\cT(f)$ we assign the number
\begin{equation}\label{multq}
m_f(B)=\nu(B)(t_f^{(1)}(B)+n(B)t_f^{(2)}(B)-1)\;
\end{equation}
Let us observe that $m_f(B)$ is positive if and only if $B\in\cE(f)$.

Eggers~\cite{E,GB} proved that
\begin{equation}\label{Qf}
Q(f)=\{q_f(B):\,B\in\cE(f)\}\;.
\end{equation}
Because of a difference in approach the analogous formulas of Eggers have a
different form. We reprove the result of Eggers in~Corollary~\ref{cor_Eggers}.
Eggers also obtained the ``multiplicities'' of $q\in Q(f)$ as the sum of
$m_f(B)$ over balls $B\in\cE(f)$ that lead to the same value of $q$
(see: Remark~\ref{rem_mult}).

Returning to Example~\ref{ex26} we obtain
$O_f(\cB)=\ord1\,f=5$, $O_f(B_1)=2$, $O_f(B_2)=3$,
$q_f(\cB)=O_f(\cB)d(\cB)=5$, $q_f(B_1)=q_f(\cB)+(d(B_1)-d(\cB))O_f(B_1)=8$,
$q_f(B_2)=q(\cB)+(d(B_2)-d(\cB))O_f(B_2)=8$, $m_f(\cB)=1$,
$m_f(B_1)=1$, $m_f(B_2)=2$.

We can apply the formulas of $Q(f)$ to compute the \L{}ojasiewicz exponent.
From (\ref{Te}) we have
\begin{equation}\label{eq_teis_eggers}
\cL_0(f)=\begin{array}[t]{c}\max\\{\scriptstyle B\in\cE(f)}\end{array}%
(q_f(B)-1)\;.
\end{equation}
In Example~\ref{ex26} we obtain $\cL_0(f)=7$ (as earlier).
Now, we can prove part (i) of the Main Result.

\vspace{1ex}\noindent{\bf Proof of Theorem~\ref{main_result}~(i)}\\
For any unitangent $f$ the formula is obvious. Assume that $f$
is multitangent. Then $\cB\in\cE(f)$. If $f=0$ is an ordinary
singularity then $\cE(f)=\{\cB\}$ and the formula is strightforward.
Let $f^{(1)},\dots,f^{(s)}$ be all the
singular unitangent components ($s\geq 1$). We have
$\cE(f)=\{\cB\}\cup\cE(f^{(1)})\cup\dots\cup\cE(f^{(s)})$ by
Property~\ref{proper_tangential}. Let us observe that
$q_f(\cB)=\ord1\,f$ and $q_f(B)>\ord1\,f$ for $B\in\cE(f)\setminus\{\cB\}$
by (\ref{monotonicity1}). Therefore, we can omit $\cB$ in
(\ref{eq_teis_eggers}). Observing that
$q_f(B)=q_{f^{(i)}}(B)+\ord1\,f-\ord1\,f^{(i)}$ for every $B\in\cE(f^{(i)})$,
$i=1,\dots,s$, we obtain
\begin{eqnarray}
\cL_0(f) & = &
\begin{array}{c}{\scriptstyle s}\\\max\\{\scriptstyle i=1}\end{array}
\begin{array}[t]{c}\max\\{\scriptstyle B\in\cE(f^{(i)})}\end{array}
(q_{f^{(i)}}(B)-1+\ord1\,f-\ord1\,f^{(i)})\\
& = &
\begin{array}{c}{\scriptstyle s}\\\max\\{\scriptstyle i=1}\end{array}
(\cL_0(f^{(i)})+\ord1\,f-\ord1\,f^{(i)})
\end{eqnarray}
and clearly we can take the last maximum over $i=1,\dots,t$ \rule{1ex}{1ex}

\section{Polar quotients and \L{}ojasiewicz exponent}\label{Loj_rel}

Let us consider a singularity $f=0$ and a regular parameter $\lambda$.
In this section we give a formula for the maximal polar quotient
$q_0(f,\lambda)$ (Corollary~\ref{cor_q}) and the formula for the
\L{}ojasiewicz exponent $\cL_0(f|\Gamma_{f,\lambda})$
(Corollary~\ref{cor_Loj_rel}). In both formulas we use the Eggers collection
extended by the balls encoding the positions of branches of the polar $\bJ(\lambda,f)=0$
with respect to the singularity $f=0$. We give three examples concerning
the inequalities from Theorem~\ref{thm_Loj_and_polar}~(a).
The most important is Example~\ref{ex_Loj_irr}. It shows that
there exists a specific equisingularity class of the pair $(f,\lambda)$
such that the \L{}ojasiewicz exponent $\cL_0(f|\Gamma_{f,\lambda})$
with respect to the polar curve
is not an equisingularity invariant inside this class
(Remark~\ref{rem_exception}).

\vspace{1ex}\noindent{\bf Position of a branch with respect to a germ}\\
We need to describe a position of a branch $h\in\cB$ with
respect to a germ $f=0$ by using equisingularity information of the
pair $(f,h)$. To this end we consider the chain
$\oK_f(h)=\{B(f_1,h),\dots,B(f_r,h)\}$ (we write $K_f(h)$ when we omit
the ball with infinite diameter). Let us denote $B_f(h)=\max\,\oK_f(h)$.
We have $d(B_f(h))<\infty$ if and only if $h$ is not a branch of $f$.
By using ($D_3$) we obtain
\begin{proper}\label{proper_extension}
$\oE(f)\cup\oK_f(h)=\oE(f)\cup\{B_f(h)\}$.
\end{proper}
\begin{proper}\label{qb}
If $h$ is not a branch of $f$ then
$$
\frac{(f,h)_0}{\ord1\,h}=q_f(B_f(h))\;.
$$
\end{proper}

\vspace{1ex}\noindent{\bf Polar quotients}\\
Now, let us consider the factorizations
\begin{equation}\label{factorizations}
f=\lambda^\delta\tilde f\quad\mbox{ and }\quad
\bJ(\lambda,f)=\lambda^{\delta}g_1\dots g_u
\end{equation}
where $\lambda$ does not divide $\tilde f$ and $g_1,\dots,g_u$ are
irreducible factors of $\bJ(\lambda,f)$ different from $\lambda$.
It is important that in both formulas we have the same
$\delta$. We will denote this number by $\delta_\lambda(f)$.
Since $f$ is reduced $\delta_\lambda(f)\in\{0,1\}$.
\begin{rem}\label{Q_empty}{\rm
Since
\begin{equation}\label{sum_mult}
\sum_{j=1}^u(g_j,\lambda)_0=(\tilde f,\lambda)_0-1\;,
\end{equation}
the condition $Q(f,\lambda)=\varnothing$ is equivalent to
$(\tilde f,\lambda)_0=1$. This means that $\tilde f$ is a smooth branch which
is transverse to $\lambda$. Hence $f=0$ is a Morse singularity with $\lambda$
as a branch.
} 
\end{rem}
\begin{cor}\label{cor_q}
\mbox{}
\begin{itemize}
\item[\rm(a)] $\displaystyle Q(f,\lambda)=\left\{
\frac{q_f(B_f(g_j))}{d(g_j,\lambda)}:\,j=1,\dots,u\right\}\;.$
\item[\rm(b)]
If $(\tilde f,\lambda)_0>1$ then
$\displaystyle
  q_0(f,\lambda)=
  \begin{array}[t]{c}\max\\
  {\scriptstyle j=1,\dots,u}
  \end{array}\frac{q_f(B_f(g_j))}{d(g_j,\lambda)}\;.
$
\end{itemize}
\end{cor}
Proof. We apply Property~\ref{qb} to (\ref{quotients}) and we use
$$
\frac{(f,g_j)_0}{(\lambda,g_j)_0}=
\frac{(f,g_j)_0}{\ord1\,g_j}\cdot
\frac{\ord1\,g_j}{(\lambda,g_j)_0}=
\frac{q_f(B_f(g_j))}{d(g_j,\lambda)},\quad\quad j=1,\dots,u\;.
$$
In Corollary~\ref{cor_Loj_rel} we obtain analogous
formulas for $\cL_0(f|\Gamma_{f,\lambda})$.

\vspace{1ex}\noindent
{\bf The \L{}ojasiewicz exponent with respect to the polar curve}\\
In the following proposition we use a natural extension of the intersection
multiplicity to quotients of series.
\begin{prop}\label{prop_Loj_rel}
Let us consider an isolated singularity $f=0$
and a regular parameter $\lambda$. Then
$$
  \cL_0(f|\Gamma_{f,\lambda})=\begin{array}[t]{c}
  \max\\{\scriptstyle h}\end{array}
  \frac{\left({\frac f\lambda},h\right)_0}{\ord1\,h}
$$
where $h$ runs over irreducible factors of the polar $\bJ(\lambda,f)$.
\end{prop}
Proof. We apply the formula from~\cite{LT}:
\begin{equation}\label{eq_LT}
  \cL_0(f|\Gamma_{f,\lambda})=
  \begin{array}[t]{c}\max\\{\scriptstyle\gamma}\end{array}
  \frac{\ord1\,((\grad\,f)\circ\gamma)}{\ord1\,\gamma}
\end{equation}
where $\gamma(T)\in\bbC\{T\}^2$, $\gamma(0)=0\in\bbC^2$, runs over a finite set of
analytic arcs that parametrize the branches of $\Gamma_{f,\lambda}$.
We can assume that $\lambda=X$. Then $\bJ(\lambda,f)=\partial f/\partial Y$.
By~(\ref{eq_LT}) we have
$$
\cL_0(f|\{\P f/\P Y=0\}) = \begin{array}[t]{c}
\max\\{\scriptstyle h}\end{array}\,
\min\left\{
\frac{\left(\frac{\P f}{\P X},h\right)_0}{\ord1\,h},
\frac{\left(\frac{\P f}{\P Y},h\right)_0}{\ord1\,h}
\right\}
=\begin{array}[t]{c}
\max\\{\scriptstyle h}\end{array}
\frac{\left(\frac{\P f}{\P X},h\right)_0}{\ord1\,h}\;
$$
where $h$ runs over all branches of $\P f/\P Y$.
Let us write $f=X^\delta\tilde f$ where $\delta=\delta_X(f)$.
If $\delta=1$ then $h=X$ is a branch of $\P f/\P Y$. In this case
$$
  \frac{\left(\frac{\P f}{\P X},h\right)_0}{\ord1\,h}=
  \frac{\left(\tilde f+X\frac{\P\tilde f}{\P X},X\right)_0}{\ord1\,X}=
  \frac{\left(\frac{f}{X},h\right)_0}{\ord1\,h}\;.
$$
When $h\neq X$ we finish by using a parametrization
of $h$ of the type $\gamma(T)=(T^N,z(T))\in\bbC\{T\}^2$,
$\gamma(0)=0$ (see~\cite{P2})~$\rule{1ex}{1ex}$

\vspace{1ex}When $\lambda$ divides $f$ we define
$\tilde B=\max\{B(\lambda,f_i):\,i=1,\dots,r,\,f_i\neq\lambda\}$.
\begin{cor}\label{cor_Loj_rel}
With notation from (\ref{factorizations})
the number $\cL_0(f|\Gamma_{f,\lambda})$ equals
{\small
\begin{itemize}
\item[\rm(1)] $\displaystyle
  \begin{array}[t]{c}\max\\
  {\scriptstyle j=1,\dots,u}
  \end{array}
  \left(q_f(B_f(g_j))-d(g_j,\lambda)\right)$ if $\delta_\lambda(f)=0$,
\item[\rm(2)]
$\displaystyle
  \max\left\{q_f(\tilde B)-d(\tilde B),\,
  \begin{array}[t]{c}\max\\
  {\scriptstyle j=1,\dots,u}
  \end{array}
  \left(q_f(B_f(g_j))-d(g_j,\lambda)\right)\right\}
$ if $\delta_\lambda(f)=1$ and $(\tilde f,\lambda)_0>1$,
\item[\rm(3)] $q_f(\tilde B)-d(\tilde B)$ if $\delta_\lambda(f)=1$
and $(\tilde f,\lambda)_0=1$.
\end{itemize}
} 
\end{cor}
Proof. By Proposition~\ref{prop_Loj_rel} and Property~\ref{qb} for $h=g_j$,
$j\in\{1,\dots,u\}$ we obtain
$$
\frac{\left(\frac{f}{\lambda},g_j\right)_0}{\ord1\,g_j}
=\frac{(f,g_j)_0}{\ord1\,g_j}-\frac{(\lambda,g_j)_0}{\ord1\,g_j}
=q_f(B_f(g_j))-d(g_j,\lambda)\;.
$$
If $h=\lambda$ then
$$
\frac{\left(\frac{f}{\lambda},\lambda\right)_0}{\ord1\,\lambda}
=(\tilde f,\lambda)_0
=q_f(\tilde B)-d(\tilde B)\;\rule{1ex}{1ex}
$$

\vspace{1ex}\noindent{\bf Examples}\\
In the first example we illustrate formulas from Corollaries~\ref{cor_q}
and~\ref{cor_Loj_rel}. We consider irreducible factors of the type
$$
aX^p+bY^q+\sum_{\alpha q+\beta p>pq}c_{\alpha\beta}X^\alpha Y^\beta,\quad
ab\neq 0,\quad\GCD(p,q)=1\;.
$$
We write shorter $aX^p+bY^q+\dots$.
\begin{ex}{\rm
Let
$$
f=Y^7+XY^4+X^2Y^2-2X^3=(Y^2-X+\dots)(Y^2+2X+\dots)(Y^3+X+\dots)
$$
and $\lambda=X$. Then
$$
\bJ(\lambda,f)=\frac{\P f}{\P Y}
=7Y^6+4XY^3+2X^2Y=Y(Y^2+\frac 12 X+\dots)(7Y^3+4X+\dots)\;.
$$
Let us denote by $f_1,f_2,f_3$ the branches of $f=0$ and by $g_1,g_2,g_3$
the branches of $\bJ(\lambda,f)=0$, respectively.
The collection $\cE(f)$ has the only one ball
$B_1=B(f_1,f_2)=B(f_1,f_3)=B(f_2,f_3)$. From~(\ref{eq_teis_eggers}) we have
$\cL_0(f)=q_f(B_1)-1=5$.
We consider the extended
collection
$$
\oE(f)\cup\{B_f(\lambda),B_f(g_1),B_f(g_2),B_f(g_3)\}
$$
and its graphical representation.
$$
\epsfbox{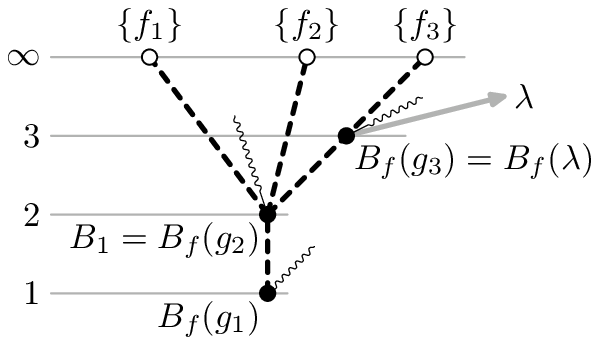}
$$
We denote the position of $B_f(\lambda)$ by an arrow and the positions
of $B_f(g_j)$ by coils. We have
$q_f(B_f(g_1))=3$,
$q_f(B_f(g_2))=6$,
$q_f(B_f(g_3))=7$,
$d(g_1,\lambda)=1$,
$d(g_2,\lambda)=2$,
$d(g_3,\lambda)=3$.
By Corollary~\ref{cor_Loj_rel} $\cL_0(f|\Gamma_{f,\lambda})=4$ and
by Corollary~\ref{cor_q} $q_0(f,\lambda)=3$.
Both inequalities from~Theorem~\ref{thm_Loj_and_polar}~(a) are strict.
} %
\end{ex}
In the following example $\lambda$ is a branch of $f$.
\begin{ex}\label{lambda}{\rm
Let us consider $f=f_1f_2=X(Y^2+X)$, $\lambda=X$.
We have $\cE(f)=\{\tilde B\}$ where $\tilde B=B(f_1,f_2)$,
$q_f(\tilde B)=4$. Hence $\cL_0(f)=q_f(\tilde B)-1=3$.
We have $\bJ(\lambda,f)=\P f/\P Y=2XY=\lambda g$.
Since $\lambda$ is a branch of $\bJ(\lambda,f)$ we denote it by a coil arrow.
$$
\epsfbox{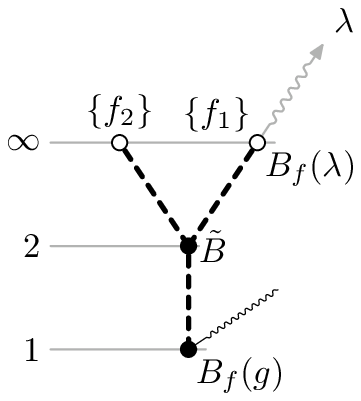}
$$
We have
$q_f(B_f(g))=2$,
$d(g,\lambda)=1$,
$d(\tilde B)=2$.
By Corollary~\ref{cor_Loj_rel}~(2)
$\cL_0(f|\Gamma_{f,\lambda})=\max\{q_f(\tilde B)-d(\tilde B),\,%
q_f(B_f(g))-d(g,\lambda)\}=\max\{2,1\}=2$ and
by Corollary~\ref{cor_q} $q_0(f,\lambda)=2$. Let us notice that here
$\cL_0(f|\Gamma_{f,\lambda})=\cL_0(f|\{\lambda=0\})$.
} 
\end{ex}
The following example shows that the position of $B_f(g_j)$
is not determined (in general) by the equisingularity class of $(f,\lambda)$
(compare~\cite{KP2}, Example~8.1). This phenomenon enable us to find
equisingular pairs $(f,\lambda)$, $(f',\lambda')$ such that
$\cL_0(f|\Gamma_{f,\lambda})\neq\cL_0(f'|\Gamma_{f',\lambda'})$.
\begin{ex}\label{ex_Loj_irr}{\rm
Let us consider $f=f_1f_2=Y^4-X^2$ and $f'=f_1'f_2'=Y^4-X^2+X^2Y$. We put
$\lambda=\lambda'=X$. We have
$\bJ(\lambda,f)=\P f/\P Y=4Y^3$ and
$\bJ(\lambda',f')=\P f'/\P Y=4Y^3+X^2$. By Corollary~\ref{cor_Loj_rel}
$\cL_0(f|\Gamma_{f,\lambda})=1$ wherease
$\cL_0(f'|\Gamma_{f',\lambda'})=\frac 32$.
$$
\epsfbox{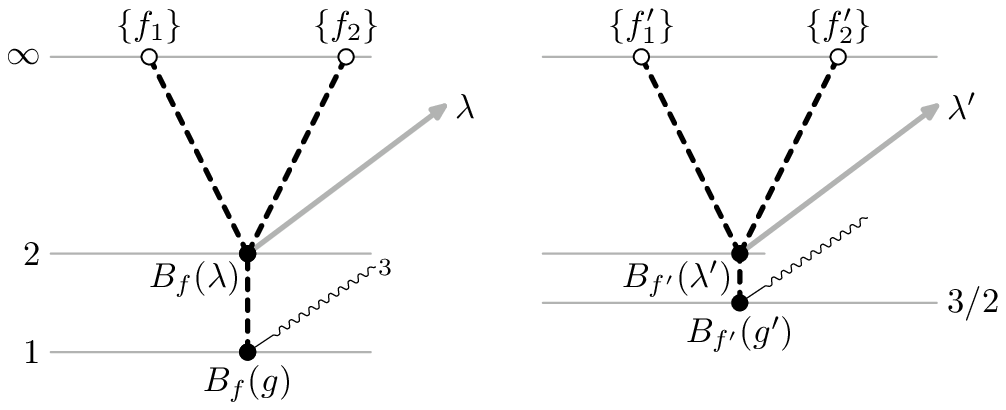}
$$
We have $\cL_0(f)=\cL_0(f')=3$ and $q_0(f,\lambda)=q_0(f',\lambda')=2$.
} %
\end{ex}
\begin{rem}\label{rem_exception}{\rm
The equisingularity class in the above example is very specific.
For the pair $(f,\lambda)$ it can be written as $t(f)=1$ and
$\cE(f)=\{B_f(\lambda)\}$. As we will see in Lemma~\ref{lm_Loj_rel} for every
different equisingularity class the \L{}ojasiewicz exponent
$\cL_0(f|\Gamma_{f,\lambda})$ is an invariant.
}
\end{rem}

\section{Factorization of polar curve}\label{factorization}

We consider a singular germ $f=0$, a regular parameter $\lambda$
and a factorization
$\bJ(\lambda,f)=\lambda^\delta g_1\dots g_u$, $\delta=\delta_\lambda(f)$,
as in (\ref{factorizations}). In this section we present 
Theorem~\ref{thm_factorization} in which every $g_j$ ($j=1,\dots,u$) is
assigned to a ball $B\in\cE(f)\cup\{B_f(\lambda)\}$ of finite diameter%
\footnote{$d(B_f(\lambda))=\infty\Leftrightarrow\delta_\lambda(f)>0$}.
This assignement corresponds to a partition $\{1,\dots,u\}=\bigcup_BJ_B$.
By putting $h_B=\prod_{j\in J_B}g_j$ we result in the factorization
$\bJ(\lambda,f)=\lambda^\delta\prod_{B}h_B$.  For $\lambda$ transversal to $f$
we obtain a version of the result of Eggers (Corollary~\ref{cor_Eggers}). Then
we describe $\cL_0(f|\Gamma_{f,\lambda})$ in terms of
the equisingularity class of pair $f,\lambda$ (Lemma~\ref{lm_Loj_rel}) and we
prove parts (ii,iii) of Main Result~A (Theorem~\ref{main_result}). Next, we
compute the polar quotients $Q(f,\lambda)$ (Proposition~\ref{prop_quotients})
and their multiplicities (Remark~\ref{rem_mult}). We describe the maximal
polar quotient $q_0(f,\lambda)$ (Lemma~\ref{polar_and_tree}). By using
Lemmas~\ref{lm_Loj_rel} and~\ref{polar_and_tree} we prove Main Result~B
(Theorem~\ref{thm_Loj_and_polar}).

\vspace{1ex}\noindent
{\bf Contact of two branches with respect to a germ}\\
Let $f_1,\dots,f_r\in\cB$ be branches of the germ $f=0$. For any $g,h\in\cB$
by ($D_3$) we obtain
\begin{equation}\label{contact_rel}
d(g,h)\geq
\begin{array}[t]{c}\max\\{\scriptstyle i=1,\dots,r}\end{array}\,
\min\{d(f_i,g),\,d(f_i,h)\}\;.
\end{equation}
We say that the contact between branches
$g$ and $h$ is {\it determined by their positions with respect to\/} $f=0$
when we have the equality in (\ref{contact_rel}). We denote the right side
of (\ref{contact_rel}) by $d_f(g,h)$.

For a ball $B\subset\cB$ and a branch $h\in\cB$ we define
\begin{equation}\label{eq_dBh}
d(B,h)=\inf\{d(h',h):\,h'\in B\}\;.
\end{equation}
\begin{proper}\label{dBh} If $B=B(g,R)$ then
$d(B,h)=\min\{d(g,h),\,R\}$.
\end{proper}
By using ($D_3$) we obtain
\begin{prop}\label{dfgh}
$d_f(g,h)=d(B_f(g),h)=d(B_f(h),g)$.
\end{prop}
{\bf Factorization theorem}\\
The role of balls from the chain
$\oK_f(\lambda)$ is specific. In order to recognize the minimal and the
maximal ball in this chain we define characteristic functions
$$
\sigma^{\mins}_{f,\lambda}(B)=\left\{\begin{array}{ll}
1 &\mbox{ if }B=\min\,\oK_f(\lambda)\\
0 &\mbox{ otherwise}\rule{0ex}{2ex}\\
\end{array}\right.
\quad\quad
\sigma^{\maxs}_{f,\lambda}(B)=\left\{\begin{array}{ll}
1 &\mbox{ if }B=B_f(\lambda)\\
0 &\mbox{ otherwise}\rule{0ex}{2ex}\\
\end{array}\right.\;.
$$
\begin{thm}\label{thm_factorization}
Let $f=0$ be a singular germ and let $\lambda$ be a regular parameter.
There exists a factorization
$\bJ(\lambda,f)=\lambda^\delta\prod_Bh_B$
where $\delta=\delta_\lambda(f)$ and $B$ runs over all the balls from
$\cE(f)\cup\{B_f(\lambda)\}$ with finite diameters such that
\begin{itemize}

\item[\rm(i)] if $g$ is a branch
of $h_B$ then $B_f(g)=B$ for $B\neq\min\,\oK_f(\lambda)$ and
$B_f(g)\leq B$ for $B=\min\,\oK_f(\lambda)$,

\item[\rm(ii)] $d(g,\lambda)=d_f(g,\lambda)$,

\item[\rm(iii)] the number $(h_B,\lambda)_0$ equals:
\begin{itemize}

\item[1.]
$d(B,\lambda)\nu(B)[t_f^{(1)}(B)+n(B)t_f^{(2)}(B)-1]\rule{0ex}{3ex}$
for $B\in\cE(f)\setminus\oK_f(\lambda)$;

\item[2.]
$d(B)n(B)[t_f(B)-1+\sigma^{\maxs}_{f,\lambda}(B)]-\sigma^{\mins}_{f,\lambda}(B)\rule{0ex}{3ex}$
for $B\in K_f(\lambda)$.

\end{itemize}
\end{itemize}
\end{thm}
We prove this theorem in Section~\ref{factorization_proof}.
We denote the number described in (iii) by $m_{f,\lambda}(B)$.
From part (iii) of Theorem~\ref{thm_factorization} we obtain
\begin{cor}\label{zero_case}
Let $B\in\cE(f)\cup\{B_f(\lambda)\}$, $d(B)<\infty$. Then $m_{f,\lambda}(B)=0$
if and only if $d(B)=1$ and one of the following conditions holds
\begin{itemize}
\item[\rm(a)] $\lambda$ is transversal to $f$ and $t(f)=1$,
\item[\rm(b)] $\lambda$ is tangent to $f$ and $t(f)=2$.
\end{itemize}
\end{cor}
For $\lambda$ transversal to $f$ we obtain
\begin{cor}\label{cor_Eggers}
Let $f=0$ be a singular germ and let $\lambda$ be a regular parameter transversal
to $f$. Then there exists a factorization $\bJ(\lambda,f)=\prod_{B\in\cE(f)}h_B$,
such that if $g$ is a branch of $h_B$ then $d(g,\lambda)=1$ and
$(f,g)_0/(\ord1\,g)=q_f(B)$. Moreover
\begin{equation}\label{Eggers_mult}
\ord1\,h_B=\nu(B)[t_f^{(1)}(B)+n(B)t_f^{(2)}(B)-1]\;.
\end{equation}
\end{cor}
Proof. We have $\min\,\oK_f(\lambda)=B_f(\lambda)=\cB$ and $\delta_\lambda(f)=0$.
The factorization from Theorem~\ref{thm_factorization} has the form
$\bJ(\lambda,f)=\prod_{B\in\cE(f)}h_B$. We can omit the ball $B_f(\lambda)$:
if $t(f)=1$ then $m_{f,\lambda}(B_f(\lambda))=0$ by Corollary~\ref{zero_case};
if $t(f)>1$ then $B_f(\lambda)\in\cE(f)$. Let $B\in\cE(f)$ and let $g$ be an
irreducible factor of $h_B$. Since the ball $\min\,\oK_f(\lambda)$ has the
minimal possible diameter, we have $B_f(g)=B$. Let us write $B=B(f_i,R)$,
$i\in\{1,\dots,r\}$. By Property~\ref{dBh} we have
$d(B,\lambda)=\min\{d(f_i,\lambda),\,R\}=1$.
From Theorem~\ref{thm_factorization}~(ii) and Proposition~\ref{dfgh} we obtain
$d(g,\lambda)=d_f(g,\lambda)=d(B_f(g),\lambda)=d(B,\lambda)=1$. Hence
$(h_B,\lambda)=\ord1\,h_B$. The equality $(f,g)_0/(\ord1\,g)=q_f(B)$ follows
directly from Property~\ref{qb}. To check (\ref{Eggers_mult}) let us observe
that $\oK_f(\lambda)=\{B_f(\lambda)\}=\{\cB\}$. Therefore, if $B\in\oK_f(\lambda)$
then $\sigma_{f,\lambda}^{\maxs}(B)=\sigma_{f,\lambda}^{\mins}(B)=1$,
$d(B,\lambda)=d(B)=1$ and $\nu(B)=n(B)=1$. In this case both formulas from
Theorem~\ref{thm_factorization}~(iii) coincide~\rule{1ex}{1ex}

\vspace{1ex}\noindent{\bf Consequences for the \L{}ojasiewicz exponent}\\
Let us observe that
for $B_1,B_2\in\cE(f)\cup\{B_f(\lambda)\}$ of finite diameters
such that $B_1<B_2$ we have
\begin{equation}\label{monotonicity}
q_f(B_1)-d(B_1,\lambda)\leq q_f(B_2)-d(B_2,\lambda)\;.
\end{equation}
\begin{lm}\label{lm_Loj_rel}
Let $f=0$ be a singular germ and let $\lambda$ be a regular parameter.
Then
$$
\cL_0(f|\Gamma_{f,\lambda})\leq%
\begin{array}[t]{c}\max\\{\scriptstyle B\in\cE(f)}\end{array}%
(q_f(B)-d(B,\lambda))\;.
$$
If $t(f)\neq 1$ or $\cE(f)\neq\{B_f(\lambda)\}$ then the equality holds.
\end{lm}
Proof. Let us assume first that $\lambda$ is not a branch of $f$.
Let us denote by $L_1$ the number from
Corollary~\ref{cor_Loj_rel}~(1). We want to show
\begin{equation}\label{L1}
  L_1\leq\begin{array}[t]{c}\max\\{\scriptstyle B\in\cE(f)}\end{array}%
(q_f(B)-d(B,\lambda))\;.
\end{equation}
with equality when $t(f)\neq 1$ or $\cE(f)\neq\{B_f(\lambda)\}$.
In order to prove (\ref{L1}) let us choose a branch
$g_j$ of $\bJ(\lambda,f)$ as in the beginning of this section.
It suffices to find a ball $B'\in\cE(f)$ such that
\begin{equation}\label{leq1}
q_f(B_f(g_j))-d(g_j,\lambda)\leq q_f(B')-d(B',\lambda)\;.
\end{equation}
Let us choose $B\in\cE(f)\cup\{B_f(\lambda)\}$ such that $g_j$
is a factor of $h_B$ from Theorem~\ref{thm_factorization}.
If $B\in\cE(f)$ then we put $B'=B$ and we show (\ref{leq1})
by using parts (i),(ii) of Theorem~\ref{thm_factorization},
Proposition~\ref{dfgh} and (\ref{monotonicity}).
If $B=B_f(\lambda)\notin\cE(f)$ then we define
$\cZ^*=\{B\in\cE(f):\,B_f(\lambda)<B\}$. When $\cZ^*$ is nonempty we choose
$B'\in\cZ^*$ and we obtain (\ref{leq1}) as earlier. When $\cZ^*=\varnothing$
we define $\cZ_*=\{B\in\cE(f):B<\,B_f(\lambda)\}$. Since $\cE(f)\neq\varnothing$
therefore $\cZ_*\neq\varnothing$. We put $B'=\max\,\cZ_*$.
In this case $B\neq\min\,\oK_f(\lambda)$. Therefore $B_f(g_j)=B$.
$$
\epsfbox{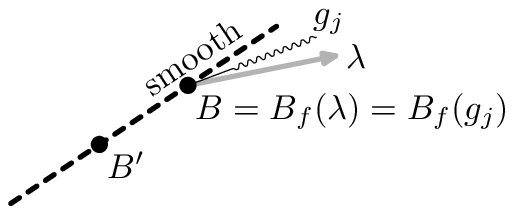}\quad
$$
Since $\char(B)\subset\char(\lambda)=\varnothing$ and $\cZ^*=\varnothing$
we have $O_f(B')=O_f(B)=1$. Therefore, in this case
we even obtain the equality in (\ref{leq1}).

In order to prove ``$\geq$'' in (\ref{L1}) let us assume that
$t(f)\neq 1$ or $\cE(f)\neq\{B_f(\lambda)\}$.
For any $B\in\cE(f)$ it suffices to
find a branch $g_j$ of $J(\lambda,f)$ such that
\begin{equation}\label{geq1}
q_f(B_f(g_j))-d(g_j,\lambda)\geq q_f(B)-d(B,\lambda)\;.
\end{equation}
If $B\neq\min\,\oK_f(\lambda)$ then $d(B)>1$ and by
Corollary~\ref{zero_case} we have $(h_B,\lambda)_0>0$.
We choose a branch $g_j$ of $h_B$ and we obtain (\ref{geq1})
(even equality) as earlier. Assume that $B=\min\,\oK_f(\lambda)$.
If $t(f)>2$ or $\lambda$ is transversal to $f$ then
$B=\min\,\oK_f(\lambda)=\cB$. Since $(h_B,\lambda)_0>0$ in this case,
we choose a branch $g_j$ of $h_B$. We have $B_f(g_j)=B$.
When $t(f)\leq 2$ and $\lambda$ is tangent to $f$ we consider two
cases. If $\min\,\oK_f(\lambda)<B_f(\lambda)$ then we choose $g_j$
as a branch of $h_{B_f(\lambda)}$. If $\min\,\oK_f(\lambda)=B_f(\lambda)$
then $f$ must be unitangent and therefore $\cE(f)\neq\{B_f(\lambda)\}$.
We take a ball $B'\in\cE(f)\setminus\{B_f(\lambda)\}$ and we choose $g_j$
as a branch of $h_{B'}$. In all these cases we obtain~(\ref{geq1}).
Hence we showed equality in~(\ref{L1}).

Now, let us assume that $\lambda$ is a branch of $f$.
Let $L_2$ be the number from Corollary~\ref{cor_Loj_rel}~(2). We want to show
\begin{equation}\label{L2}
  L_2\leq\begin{array}[t]{c}\max\\{\scriptstyle B\in\cE(f)}\end{array}%
(q_f(B)-d(B,\lambda))\;.
\end{equation}
with equality when $t(f)\neq 1$ or $\cE(f)\neq\{B_f(\lambda)\}$. We have
$\tilde B\in\cE(f)$. Then the term $q_f(\tilde B)-d(\tilde B)$ of $L_2$
is less than or equal to the right side of (\ref{L2}).
Now, let us consider a branch $g_j$ of $\bJ(\lambda,f)$
different from $\lambda$.
By Theorem~\ref{thm_factorization}
we choose $B\in\cE(f)$ such that $g_j$ is a branch of $h_B$
(we omit $B_f(\lambda)$ because $d(B_f(\lambda))=\infty$).
This $B$ gives us the expected estimation.

Now let us assume that $t(f)\neq 1$ or $\cE(f)\neq\{B_f(\lambda)\}$.
In order to prove ``$\geq$'' in (\ref{L2}) we choose $B\in\cE(f)$.
When $B=\min\,\oK_f(\lambda)$ then $B\leq\tilde B$, hence
$q_f(\tilde B)-d(\tilde B)\geq q_f(B)-d(B,\lambda)$.
If $B\neq\min\,\oK_f(\lambda)$ we consider as previously~\rule{1ex}{1ex}
\begin{rem}{\rm
Let us denote by $\ccL_0(f,\lambda)$ the number that stands on the right
side of the inequality in Lemma~\ref{lm_Loj_rel}. Clearly, it is
an equisingularity invariant of the pair $(f,\lambda)$. 
By Lemma~\ref{lm_Loj_rel} we have
$\cL_0(f|\Gamma_{f,\lambda})\leq\ccL_0(f,\lambda)$ with equality
when $t(f)\neq 1$ or $\cE(f)\neq\{B_f(\lambda)\}$.
In Example~\ref{ex_Loj_irr} we have $\ccL_0(f,\lambda)=2$.
} %
\end{rem}
\vspace{1ex}\noindent
By using Lemma~\ref{lm_Loj_rel} we can finish the proof of the Main Result.
First, we prove the following
\begin{prop}$\!\!.$\label{prop_unitangent}
Let $f=0$ be a singular unitangent germ and let $\lambda$ be a regular
parameter. Then
\begin{itemize}
\item[\rm(a)] if $\lambda$ is tangent to $f$ then
$\cL_0(f|\Gamma_{f,\lambda})\leq\ccL_0(f,\lambda)<\cL_0(f)$,
\item[\rm(b)] if $\lambda$ is transversal to $f$ then
$\cL_0(f|\Gamma_{f,\lambda})=\ccL_0(f,\lambda)=\cL_0(f)$.
\end{itemize}
\end{prop}
Proof.~(a) Since $f$ is singular we have $\cE(f)\neq\varnothing$.
Let $B_1=\min\,\cE(f)$ and let $B_2=\min\,\oK_f(\lambda)$.
We have $d(B_1)>1$ and $d(B_2)>1$. For any $B\in\cE(f)$
we obtain $d(B,\lambda)\geq\min\{d(B_1),d(B_2)\}>1$.
By Lemma~\ref{lm_Loj_rel}
$$
\cL_0(f|\Gamma_{f,\lambda})\leq\ccL_0(f,\lambda)
=\begin{array}[t]{c}\max\\{\scriptstyle B\in\cE(f)}\end{array}
(q_f(B)-d(B,\lambda))<\cL_0(f)\;.
$$
(b) If $\lambda$ is transversal to $f$ then $\cE(f)\neq\{B_f(\lambda)\}$.
Moreover $d(B,\lambda)=1$ for every $B\in\cE(f)$.
We apply Lemma~\ref{lm_Loj_rel}~\rule{1ex}{1ex}

\vspace{1ex}\noindent{\bf Proof of Theorem~\ref{main_result}~(ii),(iii)}\\
Let $f^{(1)},\dots,f^{(t)}$, $t=t(f)$, be unitangent components of $f$ and
let us denote
\begin{equation}\label{M}
M_i=\cL_0(f^{(i)})+\ord1\,f-\ord1\,f^{(i)}\;,\quad i=1,\dots,t\;.
\end{equation}
We have $M_i>\ord1\,f-1$ if and only if $f^{(i)}$ is singular ($i=1,\dots,t$).
We may assume that $M_1\geq\dots\geq M_t$. From part (i) we have
$\cL_0(f)=M_1$. Let $s$ be the number of singular components ($0\leq s\leq t$).

\vspace{1ex}\noindent Proof of (ii).
We claim that $M_1>\ord1\,f-1$. For $t(f)=1$ it follows from the fact that
$f$ is singular. For $t(f)>1$ it is a consequence of the assumption that
the maximum $M_1$ is realized for exactly one index from $\{1,\dots,t\}$.
Hence, the corresponding component $f^{(1)}$ is singular and therefore
$s\geq 1$. Let $\lambda$ be a regular parameter. As in the proof of part~(i)
of the theorem we obtain
\begin{equation}
\ccL_0(f,\lambda)=
\begin{array}[t]{c}\max\\{\scriptstyle i=1,\dots,s}\end{array}
\{\ccL_0(f^{(i)},\lambda)+\ord1\,f-\ord1\,f^{(i)}\}\;.
\end{equation}
Assume that $\lambda$ is tangent to $f^{(1)}$. If $s=1$ then by 
Proposition~\ref{prop_unitangent}~(a)
\begin{eqnarray*}
\cL_0(f|\Gamma_{f,\lambda})\leq\ccL_0(f,\lambda) & = & \ccL_0(f^{(1)},\lambda)
+\ord1\,f-\ord1\,f^{(1)}\\ & < & \cL_0(f^{(1)})
+\ord1\,f-\ord1\,f^{(1)} = \cL_0(f)\;.
\end{eqnarray*}
If $s>1$ then $M_1>M_2$ 
and we have
$$
\cL_0(f|\Gamma_{f,\lambda})\leq\ccL_0(f,\lambda) = 
\max\{\ccL_0(f^{(1)},\lambda)+\ord1\,f-\ord1\,f^{(1)},M_2\}
<\cL_0(f)\;.
$$
In order to prove the opposite implication in (ii) suppose that
$\lambda$ is transversal to $f^{(1)}$. In this case $\cB\in\cE(f)$,
therefore the condition $\cE(f)\neq\{B_f(\lambda)\}$ is satisfied.
According to Lemma~\ref{lm_Loj_rel} we have
$$
\cL_0(f|\Gamma_{f,\lambda})=\ccL_0(f,\lambda) = \cL_0(f^{(1)})
+\ord1\,f-\ord1\,f^{(1)} = \cL_0(f)\;\rule{1ex}{1ex}
$$
Proof of (iii). We have $t=t(f)\geq 2$ and $M_1=M_2=\cL_0(f)$.
If $M_1=\ord1\,f-1$, then all the tangential components of $f$
are nonsingular (ordinary singularity). In this case $\cE(f)=\{\cB\}$.
By Lemma~\ref{lm_Loj_rel}
$\cL_0(f,\lambda)=q_f(\cB)-d(\cB,\lambda)=\ord1\,f-1=\cL_0(f)$
for every regular parameter $\lambda$.
If $M_1>\ord1\,f-1$ then $f^{(1)}$ and $f^{(2)}$ are singular ($s\geq 2$).
Since every regular parameter $\lambda$ is transversal to
$f^{(1)}$ or to $f^{(2)}$ we obtain $\cL_0(f|\Gamma_{f,\lambda})=\cL_0(f)$
as earlier~\rule{1ex}{1ex}

\vspace{1ex}\noindent{\bf Consequences for polar quotients}\\
Below, we apply Theorem~\ref{thm_factorization} to polar quotients.
We use notation of this theorem.
\begin{prop}\label{prop_FAKT}
If $B=\min\,\oK_f(\lambda)$ and $g$ is a branch of $h_B$ then
$d(g,\lambda)=d(B_f(g))$.
\end{prop}
Proof. By (i) of the theorem we have $B_f(g)\leq B=\min\,\oK_f(\lambda)$. Since
$d(\min\,\oK_f(\lambda))\leq\min_{i,j}d(f_i,f_j)$, then
$d(f_1,g)=\dots=d(f_r,g)=d(B_f(g))$. By using (ii) we obtain
$d(g,\lambda)=d_f(g,\lambda)=\max_i\min\{d(f_i,g),d(f_i,\lambda)\}=d(B_f(g))$%
~\rule{1ex}{1ex}

From Corollary~\ref{cor_q} and Theorem~\ref{thm_factorization} we obtain
\begin{prop}\label{prop_quotients}
$$
Q(f,\lambda)=\left\{\frac{q_f(B)}{d(B,\lambda)}:\,
B\in\cE(f)\cup\{B_f(\lambda)\},\,d(B)<\infty,\,m_{f,\lambda}(B)>0\right\}\;.
$$
\end{prop}
Proof. Let us choose a factor $g_j$ of $\bJ(\lambda,f)$ ($j\in\{1,\dots,u\}$)
as in Corollary~\ref{cor_q}. By Theorem~\ref{thm_factorization}~(i)
$g_j$ is a factor of $h_B$ for $B\in\cE(f)\cup\{B_f(\lambda)\}$, $d(B)<\infty$,
$m_{f,\lambda}(B)>0$. If $B\neq\min\,\oK_f(\lambda)$ then $B_f(g_j)=B$.
We finish by using (ii) and Proposition~\ref{dfgh}:
$d(g_j,\lambda)=d_f(g_j,\lambda)=d(B_f(g_j),\lambda)=d(B,\lambda)$.

If $B=\min\,\oK_f(\lambda)$ then $B_f(g_j)\leq B$. We have
$q_f(B_f(g_j))=(\ord1\,f)d(B_f(g_j))$ and
$q_f(B)=(\ord1\,f)d(B)$. We finish by using Proposition~\ref{prop_FAKT}
and Property~\ref{dBh}:
$$
\frac{q_f(B_f(g_j))}{d(g_j,\lambda))}=
\frac{q_f(B_f(g_j))}{d(B_f(g_j))}=
\ord1\,f=
\frac{q_f(B)}{d(B)}=
\frac{q_f(B)}{d(B,\lambda)}\;\rule{1ex}{1ex}
$$
\begin{rem}\label{rem_mult}
{\rm(multiplicities of polar quotients)\\
To every $q\in Q(f,\lambda)$ we can assign a multiplicity
$m_q=\sum_B m_{f,\lambda}(B)$
where $B$ runs over all balls from $\cE(f)\cup\{B_f(\lambda)\}$ with finite
diameters such that $q_f(B)/d(B,\lambda)=q$.
By (\ref{sum_mult}) we have $\sum_{q\in Q(f,\lambda)}m_q=(\tilde{f},\lambda)_0-1$.
} %
\end{rem}
It is important in the following lemma that we can omit the ball
$B_f(\lambda)$ in the cases~(i) and~(ii).
\begin{lm}\label{polar_and_tree}
{\rm(description of the maximal polar quotient)}
\begin{itemize}
\item[\rm(i)] Assume that $t(f)\neq 2$. Then
\begin{equation}\label{eq1}
  q_0(f,\lambda)=%
  \begin{array}[t]{c}\max\\{\scriptstyle B\in\cE(f)}\end{array}%
  \frac{q_f(B)}{d(B,\lambda)}\;,
\end{equation}
\item[\rm(ii)] If $t(f)=2$ and $\#\cE(f)\geq 2$ then
\begin{equation}\label{eq2}
  q_0(f,\lambda)=%
  \begin{array}[t]{c}\max\\{\scriptstyle B\in\cE(f)\setminus\{\cB\}}\end{array}%
  \frac{q_f(B)}{d(B,\lambda)}\;,
\end{equation}
\item[\rm(iii)] Assume that $t(f)=2$ and $\#\cE(f)=1$ {\rm(}Morse case{\rm)}.
If $\lambda$ is not a brach of $f$ then
\begin{equation}\label{eq3}
  q_0(f,\lambda)=%
  \frac{q_f(B_f(\lambda))}{d(B_f(\lambda))}=%
  \frac{(f,\lambda)_0}{(f,\lambda)_0-1}
\end{equation}
and if $\lambda$ is a branch of $f$ then $q_0(f,\lambda)=-\infty$.
\end{itemize}
\end{lm}
Proof. Let us consider the set of balls from Proposition~\ref{prop_quotients}
$$
\cZ=\{B\in\cE(f)\cup\{B_f(\lambda)\}:\,d(B)<\infty,\,m_{f,\lambda}(B)>0\}\;.
$$
(i) In order to prove ($\leq$) in (\ref{eq1})
we choose $B\in\cZ$. It suffices to find $B'\in\cE(f)$ such that
$q_f(B)/d(B,\lambda)\leq q_f(B')/d(B',\lambda)$. If $B\in\cE(f)$ we put
$B'=B$. Suppose that $B=B_f(\lambda)$. Since $\cE(f)$ is nonempty, at least
one of the following conditions holds:
(a) there exists $B_1\in\cE(f)$ such that
$B_1<B$ and $B_1$ is the direct predecessor of $B$,
(b) there exists $B_2\in\cE(f)$
such that $B<B_2$ and $B_2$ is the direct successor of $B$.
In case (a) we have
$d(B_1,\lambda)=d(B_1)$, $d(B,\lambda)=d(B)$. Moreover, $O_f(B)=1$ because
$\lambda$ is smooth. Hence
$$
  \frac{q_f(B_1)}{d(B_1,\lambda)}-\frac{q_f(B)}{d(B,\lambda)}
  =\frac{(q_f(B_1)-d(B_1))(d(B)-d(B_1))}{d(B_1)d(B)}>0\;.
$$
In case (b) we have $d(B_2,\lambda)=d(B,\lambda)=d(B)$. Therefore,
$$
  \frac{q_f(B_2)}{d(B_2,\lambda)}-\frac{q_f(B)}{d(B,\lambda)}
  =\frac{O_f(B_2)(d(B_2)-d(B))}{d(B)}>0\;.
$$
As $B'$ we choose $B_1$ or $B_2$. In order to prove inequality ($\geq$)
it suffices to show $\cZ\supset\cE(f)$. Let $B\in\cE(f)$ and suppose that
$m_{f,\lambda}(B)=0$. Since $t(f)\neq 2$, by Corollary~\ref{zero_case} we
obtain $d(B)=1$ and $f$ is unitangent. Therefore, $d(B)>1$ for every
$B\in\cE(f)$, which is a contradiction. Hence $m_{f,\lambda}(B)>0$
and inclusion is proved.

\vspace{1ex}\noindent
(ii) Since $t(f)=2$ we have $\cB\in\cE(f)$. Let us notice that
$m_{f,\lambda}(\cB)=\sigma_{f,\lambda}^{\maxs}(\cB)$.
If $\lambda$ is tangent to $f$ then $\sigma_{f,\lambda}^{\maxs}(\cB)=0$.
Hence $\cB\notin\cZ$. By the assumption $\#\cE(f)\geq 2$ the set $\cE(f)\setminus\{\cB\}$
is nonempty. We prove ($\leq$) as in~(i).
The inequality ($\geq$) is a consequence of $\cZ\supset\cE(f)\setminus\{\cB\}$.
If $\lambda$ is transversal to $f$ then $B_f(\lambda)=\cB$ and $\cZ=\cE(f)$.
Since $\cE(f)\setminus\{\cB\}$ is nonempty,
the inequality ($\leq$) follows from the fact that $q_f(\cB)/d(\cB,\lambda)$
is now the minimal possible polar quotient. The inequality ($\geq$) follows
from the obvious inclusion as earlier.

\vspace{1ex}\noindent
(iii) We have $\cE(f)=\{\cB\}$. If $\lambda$ is transversal to $f$ then
$B_f(\lambda)=\cB$, $m_{f,\lambda}(\cB)=\sigma_{f,\lambda}^{\maxs}(\cB)=1$,
$q_f(\cB)/d(\cB,\lambda)=2$. If $\lambda$ is tangent to $f$ then
$m_{f,\lambda}(\cB)=0$. If $\lambda$ is not a branch of $f$ then
$m_{f,\lambda}(B_f(\lambda))=d(B_f(\lambda))\geq 1$. Then
$q_f(B_f(\lambda))/d(B_f(\lambda))=(f,\lambda)_0/((f,\lambda)_0-1)$. If $\lambda$
is a branch of $f$ then $d(B_f(\lambda))=\infty$ and
$Q(f,\lambda)=\emptyset$~\rule{1ex}{1ex}

\begin{cor}\label{q_all_cases} In all the cases
$\displaystyle
\begin{array}[t]{c}\max\\{\scriptstyle B\in\cE(f)}\end{array}%
\frac{q_f(B)}{d(B,\lambda)}\geq
q_0(f,\lambda)\;.
$
\end{cor}
\begin{ex} {\rm
Let $f=f_1f_2f_3=Y(Y^2-X)(Y^2+X)$, $\lambda=X$.
We have $\cE(f)=\{B_0,B_1\}$, where $B_0=B(f_1,f_2)=B(f_1,f_3)$
and $B_1=B(f_2,f_3)$ with $q_f(B_0)=3$, $q_f(B_1)=5$,
$d(B_0,\lambda)=1$ and $d(B_1,\lambda)=2$.
$$
\epsfbox{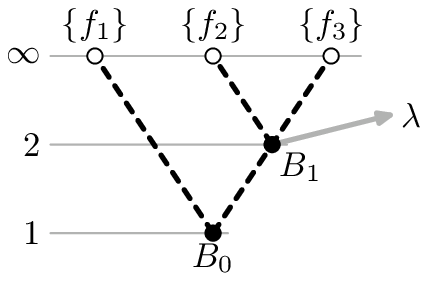}
$$
Although $q_f(B_0)/d(B_0,\lambda)>q_f(B_1)/d(B_1,\lambda)$ we have
$q_0(f,\lambda)=q_f(B_1)/d(B_1,\lambda)=5/2$ since we omit $B_0=\cB$
by Lemma~\ref{polar_and_tree}~(ii).
} 
\end{ex}

\vspace{1ex}\noindent
By using both Lemmas~\ref{lm_Loj_rel} and~\ref{polar_and_tree} we can
prove Theorem~\ref{thm_Loj_and_polar}.

\vspace{1ex}\noindent{\bf Proof of Theorem~\ref{thm_Loj_and_polar}}\\
(a) For any $B\in\cE(f)$ we have
$$
  q_f(B)-1\geq q_f(B)-d(B,\lambda)\geq\frac{q_f(B)}{d(B,\lambda)}-1\;.
$$
If $t(f)\neq 1$ or $\cE(f)\neq\{B_f(\lambda)\}$ then we finish by using
(\ref{eq_teis_eggers}), Lemma~\ref{lm_Loj_rel} and Corollary~\ref{q_all_cases}.
When $t(f)=1$ and $\cE(f)=\{B_f(\lambda)\}$ let us denote $B_0=B_f(\lambda)$.
We have $B_0=\min\,\oK_f(\lambda)$. By Theorem~\ref{thm_factorization}
$\bJ(f,\lambda)=h_{B_0}$, $(h_{B_0},\lambda)_0=d(B_0)n(B_0)t_f(B_0)-1>0$.
By using Corollary~\ref{cor_Loj_rel} and Proposition~\ref{prop_FAKT} we obtain
$$
  \cL_0(f|\Gamma_{f,\lambda})=
  \begin{array}[t]{c}\max\\{\scriptstyle g}\end{array}
  (q_f(B_f(g))-d(g,\lambda))=
  \begin{array}[t]{c}\max\\{\scriptstyle g}\end{array}
  (\ord1\,f-1)d(B_f(g))\;,
$$
where $g$ runs over irreducible factors of $h_{B_0}$.
As in the proof of Proposition~\ref{prop_quotients} we show that
$q_0(f,\lambda)=\ord1\,f$. Moreover $\cL_0(f)=(\ord1\,f)d(B_0)-1$,
which gives desired inequalities.

\vspace{1ex}\noindent (b) Let us consider three cases
{\small
\begin{itemize}
\item[(I)] $t(f)>2$ or ($t(f)=1$ and $\#\cE(f)>1$) or ($t(f)=1$ and $\cE(f)\neq\{B_f(\lambda)\}$)
\item[(II)] $t(f)=1$ and $\cE(f)=\{B_f(\lambda)\}$
\item[(III)] $t(f)=2$ and $\#\cE(f)>1$
\end{itemize}
} 
$$
\epsfbox{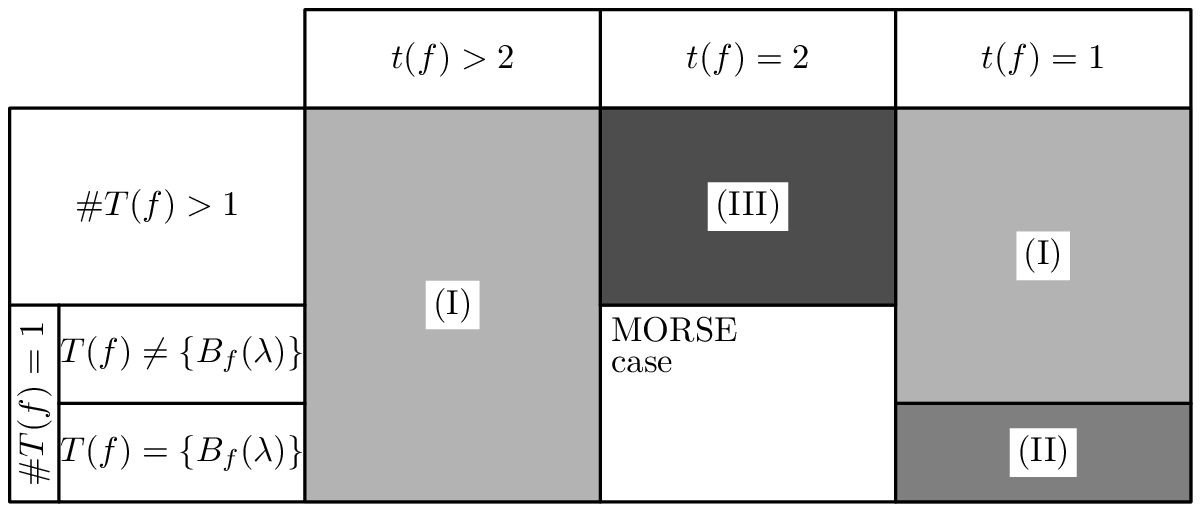}
$$
(I) In this case we obtain the desired equivalency from the formulas
$$
  \cL_0(f|\Gamma_{f,\lambda})=
  \begin{array}[t]{c}\max\\{\scriptstyle B\in\cE(f)}\end{array}
  (q_f(B)-d(g,\lambda))\;,\quad\quad
  q_0(f,\lambda)=
  \begin{array}[t]{c}\max\\{\scriptstyle B\in\cE(f)}\end{array}
  \frac{q_f(B)}{d(B,\lambda)}\;.
$$
(II) We prove in this case that none of the equalities from the statement
of the theorem can not be satisfied. Let us denote $B_0=B_f(\lambda)$.
Since the singularity is unitangent, $d(B_0)>1$.
By using the formulas as in the proof of (a) we obtain
$q_0(f)-q_0(f,\lambda)\geq(\ord1\,f)(d(B_0)-1)$
and $\cL_0(f)-\cL_0(f|\Gamma_{f,\lambda})\geq d(B_0)-1$.

\vspace{1ex}\noindent
(III) Now, the following formulas are true
$$
  \cL_0(f|\Gamma_{f,\lambda})=
  \begin{array}[t]{c}\max\\{\scriptstyle B\in\cE(f)}\end{array}
  (q_f(B)-d(g,\lambda))\;,\quad\quad
  q_0(f,\lambda)=
  \begin{array}[t]{c}\max\\{\scriptstyle B\in\cE(f)\setminus\{\cB\}}\end{array}
  \frac{q_f(B)}{d(B,\lambda)}\;.
$$
($\Rightarrow$) Assume that $\cL_0(f)=\cL_0(f|\Gamma_{f,\lambda})$. Hence,
there exists $B_1\in\cE(f)$ such that $\cL_0(f)=q_f(B_1)-d(B_1,\lambda)$,
therefore $d(B_1,\lambda)=1$. Since $\#\cE(f)>1$, and taking into consideration
(28) we can assume that $B_1\neq\cB$. For $B_1$ we obtain
$q_0(f)=q_f(B_1)=q_0(f,\lambda)$.

\vspace{1ex}\noindent($\Leftarrow$) Let us assume that $q_0(f)=q_0(f,\lambda)$.
Hence, there exists $B_2\in\cE(f)\setminus\{\cB\}$ such that
$q_f(B_2)/d(B_2,\lambda)=q_0(f)$. Therefore, $q_f(B_2)=q_0(f)$
and $d(B_2,\lambda)=1$. We obtain
$\cL_0(f)=\cL_0(f|\Gamma_{f,\lambda})$~\rule{1ex}{1ex}

\section{Proof of factorization theorem}\label{factorization_proof}

In this section we prove Theorem~\ref{thm_factorization}.
We can consider the derivative as the polar curve.
We apply the main result of~\cite{Len2004}, where a version of the Newton
algorithm~\cite{Cano} provides a description of all polar quotiens
including multiplicities (\cite{Len2004},Theorem~2.1). Here we reformulate
this result to describe the roots of the derivative (Theorem~\ref{thm_Len2004}).
Next, we study the characteristics of branches and we describe all possible
balls within the fixed coordinate system.
By using this description we assign the roots of the derivative
to the balls (Lemma~\ref{lm_balls}).
The number of these roots, described by Theorem~\ref{thm_Len2004},
gives us ``multiplicities'' of branches of the derivative assigned
to the balls (Proposition~\ref{prop_multiplicities}).

\vspace{1ex}\noindent{\bf The roots of derivative}\\
We need some preliminaries. We consider the ring
$\bbC\{X\}^*=\bigcup_{N\geq 1}\bbC\{X^{1/N}\}$ of {\it Puiseux series\/}.
For every nonzero $y(X)\in\bbC\{X\}^*$ the {\it order\/} $\ord1\,y$
stands for the minimal power with nonzero coefficient and $\In\,y$
is the corresponding monomial. We put $\ord1\,0=\infty$ and $\In\,0=0$.
It is convenient to consider the ring
$\bbC\{X^*,Y\}=\bigcup_{N\geq 1}\bbC\{X^{1/N},Y\}$.
Take $f=\sum c_{\alpha\beta}X^\alpha Y^\beta\in\bbC\{X^*,Y\}$.
As usual, we define the {\it support\/}
$\supp\,f$ as $\{(\alpha,\beta):\;c_{\alpha\beta}\neq 0\}$,
the {\it Newton diagram\/} $\Delta(f)$ as $\conv(\supp\,f+\bbR_+^2)$,
and the {\it Newton polygon\/} $\cN_f=\cN(f)$ as the set of compact faces
of $\Delta(f)$
(we use the term ``face'' in the meaning of ``1-dimensional face'').
By $\delta_Y(f)$ (resp. $\delta_X(f)$) we denote the distance
between $\Delta(f)$ and the horizontal axis (resp. vertical axis).

For $S\in\cN_f$, by $|S|_1$ and $|S|_2$ we denote the lengths of projections
of $S$ onto the horizontal and vertical axes, respectively.
We call the ratio $|S|_1/|S|_2$ the {\it inclination\/} of $S$.
We denote it by $\incl(S)$. We define
$\incl(\cN_f)=\{\incl(S):\,S\in\cN_f\}$ if $\delta_Y(f)=0$
or $\incl(\cN_f)=\{\incl(S):\,S\in\cN_f\}\cup\{\infty\}$ if $\delta_Y(f)>0$.
For $\theta>0$ (or $\theta=-\infty$) it is
useful to consider the polygon $\cN^\theta_f$ which consists of
all $S\in\cN_f$ with $\incl(S)>\theta$. We have
$\incl(\cN^\theta_f)=\incl(\cN_f)\cap(\theta,\infty]$.
We define the {\it initial form\/} of $f$
with respect to $S$ as $\In(f,S)=\sum c_{\alpha\beta}X^\alpha Y^\beta$ where
$(\alpha,\beta)\in S\cap\supp\,f$. By $t(f,S)$ we denote the number
of different roots of the polynomial $\In(f,S)(1,Y)\in\bbC[Y]$.
The number $\varepsilon(S)\in\{-1,0\}$ is defined as $-1$ when $S$ touches
the horizontal axis and as $0$, otherwise.
Put $d(f,S)=|S|_2+\varepsilon(S)-t(f,S)+1$. Note that $d(f,S)=0$
if and only if every nonzero root of $\In(f,S)$ in $\bbC\{X\}^*$
is of multiplicity $1$.
Then, we call the series $f$ {\it nondegenerate\/} on $S$.

For any $\varphi\in\bbC\{X\}^*$, $\ord1\,\varphi>0$ one can apply
the substitution $f_\varphi(X,Y)=f(X,\varphi+Y)\in\bbC\{X^*,Y\}$
(\cite{Cano}, \cite{GP1991}, \cite{KP1}).
Clearly, $f_\varphi=f$ for $\varphi=0$.
Consider the ring $\bbC[X]^*=\bigcup_{N\geq 1}\bbC[X^{1/N}]$ of
{\it Puiseux polynomials\/}.
For $\varphi\in\bbC[X]^*$, $\deg\varphi<\infty$.
Put $\deg 0=-\infty$. We define the set $\Track(f)\subset\bbC[X]^*$
of {\it tracks\/} (of the Newton algorithm) for $f$ as the minimal set
satisfying two properties: (I)~$0\in \Track(f)$, (II)~for every
$\varphi\in\Track(f)$, if there exists
$S\in\cN^{\deg\varphi}(f_\varphi)$, then for every nonzero root $aX^\theta$
of $\In(f_\varphi,S)$, $\varphi+aX^\theta\in\Track(f)$. In \cite{Len2004}
(Proposition~3.11) we give three different characterizations of the set
$\Track(f)$. We will write $\cN_\varphi$ instead of
$\cN^{\deg\varphi}(f_\varphi)$ when $f$ is fixed.
We call a series $\psi\in\bbC\{X\}^*$ a {\it continuation\/} of
$\varphi\in\bbC[X]^*$ if $\ord1\,(\varphi-\psi)>\deg\varphi$. Then we write
$\psi=\varphi+\dots\,$. By $\Track_\varphi(f)$ we denote the set of all
tracks from $\Track(f)$ that are continuations of $\varphi$.

In order to deal with multiple roots we use the notion of
{\it symmetric power\/}~\cite{Whitney}. For elements $a_1,\dots,a_s$ of a given set we
define the {\it system\/} $\bA=\langle a_1,\dots,a_s\rangle$ as the sequence
$a_1,\dots,a_s$ treated as unordered. Put $\deg\bA=s$. For
$\bA=\langle a_1,\dots,a_s\rangle$ and $\bB=\langle b_1,\dots,b_t\rangle$
we have a natural addition
$\bA\oplus\bB=\langle a_1,\dots a_s,b_1,\dots b_t\rangle$
with the neutral element $\langle\rangle$ ({\it empty system\/}).
Instead of
$\langle\underbrace{a,\dots, a}_{m\,\mbox{\scriptsize times}}\rangle$ we write
$\langle a:m\rangle$ with convention $\langle a:0\rangle=\langle\rangle$.
If $a$ appears in $\bA$ at least one time
then we write $a\in\bA$.

Now, assume that $\ord1\,f(0,Y)=p>0$. We consider
the system $\Zer\,f=\langle y_1,\dots,y_p\rangle$ of all solutions of
$f=0$ in $\bbC\{X\}^*$. Let $\varphi\in\Track(f)$. By $\Zer_\varphi{f}$
we denote the system of all solutions from $\Zer\,f$ that are continuations
of $\varphi$. Our aim is to describe the system
$\Zer(\P{f}/\P{Y})=\langle z_1,\dots,z_{p-1}\rangle$.
We define a solution $z(X)\in\Zer_\varphi(\P{f}/\P{Y})$ to be of the
$\varphi$-{\it first kind\/} if $\ord1(z-\varphi)\in\incl(\cN_\varphi)$
and of the $\varphi$-{\it second kind\/} otherwise. We control the
``kind'' by the following proposition (see: Proposition~3.4,
\cite{LenMast2009}).
For $S\in\cN_\varphi$ we put $w_{\varphi,S}(Y)=\In(f_\varphi,S)(1,Y)$.
\begin{prop}\label{prop_LenMast2009}
Let $z(X)\in\Zer_\varphi(\P{f}/\P{Y})$. Then
\begin{itemize}

\item[\rm(i)] If $z(X)$ is of the $\varphi$-first kind then:
   \begin{itemize}
   \item[\rm(a)] if $\ord1(z-\varphi)=\infty$ $($i.e. $z=\varphi)$ then $\delta_Y(f_\varphi)>1$,
   \item[\rm(b)] if $\ord1(z-\varphi)<\infty$ then there exists
   $S\in\cN(f_\varphi)$ such that $z(X)=\varphi+aX^{\incls(S)}+\dots$
   $(a\neq 0)$ and $w_{\varphi,S}'(a)=0$.
   \end{itemize}

\item[\rm(ii)] Solutions of the $\varphi$-second kind exist if and only if
both conditions hold:
  \begin{itemize}
  \item the lowest face $S=L$ of $\cN_\varphi$ touches the horizontal axis
  $($i.e. $w_{\varphi,L}(0)\neq 0)$,
  \item $\ord1(w_{\varphi,L}(Y)-w_{\varphi,L}(0))\geq 2$.
  \end{itemize}

\item[\rm(iii)] If $z(X)$ is of the $\varphi$-second kind
then $\ord1(z-\varphi)>\max\,\incl(\cN_\varphi)$.
\end{itemize}
\end{prop}
Let $\varphi\in\Track(f)$. We define the system
$\Zer_\varphi^{\fins}(\P{f}/\P{Y})$ (resp. $\Zer_\varphi^\infty(\P{f}/\P{Y})$)
which consists of those $z(X)\in\Zer_\varphi(\P{f}/\P{Y})$ that
$\ord1\,f(X,z(X))<\infty$ (resp. $\ord1\,f(X,z(X))=\infty$). We have
$\Zer_\varphi(\P{f}/\P{Y})=\Zer_\varphi^{\fins}(\P{f}/\P{Y})\oplus%
\Zer_\varphi^\infty(\P{f}/\P{Y})$.
We put $\bC_\varphi=\langle\varphi:\,\delta_Y(f_\varphi)-1\rangle$
if $\delta_Y(f_\varphi)>1$ and $\bC_\varphi=\langle\rangle$ if
$\delta_Y(f_\varphi)\in\{0,1\}$. For $S\in\cN_\varphi$ we define the system
$\bB_{\varphi,S}$ (resp. $\bA_{\varphi,S}^I$) of the $\varphi$-first kind
solutions $z(X)$ such that
$\ord1(z-\varphi)=\incl(S)$ and $\In(z-\varphi)$ is a root
(resp. is not a root) of $\In(f_\varphi,S)$. By $\bA_\varphi^{II}$ we denote
the system of all $\varphi$-second kind solutions. We put
$$
\renewcommand{\arraystretch}{1.2}
  \bA_{\varphi,S} = \left\{\begin{array}{ll}
  \bA_{\varphi,S}^I &\mbox{if $S$ does not touch the horizontal axis}\\
  \bA_{\varphi,S}^I\oplus\bA_\varphi^{II} &
  \mbox{if $S$ touches the horizontal axis}\\
  \end{array}\right.
$$
The following theorem is a reformulation of Theorem~2.1 from~\cite{Len2004}.
The proof is analogous.
\begin{thm}\label{thm_Len2004}
Let $\varphi\in\Track(f)$.
\begin{itemize}
\item[\rm(a)] $\displaystyle\Zer_{\varphi}(\P{f}/\P{Y})=
  \left[\bigoplus_{\scriptstyle S\in\cN_{\varphi}}
  (\bA_{\varphi,S}\oplus\bB_{\varphi,S})\right]\oplus\bC_{\varphi}$\\
with $\deg\bA_{\varphi,S}=t(f_\varphi,S)-1$,
$\deg\bB_{\varphi,S}=d(f_\varphi,S)$.

\item[\rm(b)] Let $S\in\cN_\varphi$; $\displaystyle\bB_{\varphi,S}=
              \bigoplus_{aX^\theta}\Zer_{\varphi+aX^\theta}(\P{f}/\P{Y})$,\\
              where $aX^\theta$ runs over all multiple nonzero roots
              of $\In(f_\varphi,S)$.
\item[\rm(c)]
$\displaystyle
  \Zer^{\fins}_\varphi(\P{f}/\P{Y})=\bigoplus_{\psi\in\Tracks_\varphi(f)}%
  \bigoplus_{\scriptstyle S\in\cN_{\psi}}
  \bA_{\psi,S}
$
\end{itemize}
\end{thm}
Now, assume that $f$ is reduced.
Then for every $\varphi\in\Track(f)$ we have
$\Zer_\varphi(\P{f}/\P{Y})=\Zer_\varphi^{\fins}(\P{f}/\P{Y})$ and
part ``$\bC$'' disappears.
For $\varphi=0$ we obtain the following two corollaries.
We write $\bA_S$, $\bB_S$ instead of $\bA_{0,S}$, $\bB_{0,S}$.
\begin{cor}\label{cor1}{\rm(see Corollary~2.5(a),\cite{Len2004})}
$$
  \Zer(\P{f}/\P{Y})=\bigoplus_{S\in\cN_f}(\bA_S\oplus\bB_S)\;.
$$
For every $S\in\cN_f$
\begin{itemize}
\item[\rm(a)] $\deg\bA_S=t(f,S)-1$,
\item[\rm(b)] $\deg\bB_S=d(f,S)$.
\end{itemize}
\end{cor}
\begin{cor}\label{cor2} For $S\in\cN_f$
$$
  \bB_S=\bigoplus_{\varphi\in\Tracks_{aX^\theta}(f)}%
  \bigoplus_{\scriptstyle S\in\cN_{\varphi}}
  \bA_{\varphi,S}\;,
$$
where $aX^\theta$ runs over all multiple nonzero roots of $\In(f_\varphi,S)$.
\end{cor}

The information presented in Corollary~\ref{cor1} corresponds to the first step
of the Newton algorithm.
The information presented in Corollary~\ref{cor2} corresponds to the
following steps.

\vspace{1ex}\noindent{\bf Characteristic Newton diagram of a branch}\\
Recall a notion of the {\it cycle\/} generated by a Puiseux series
$y(X)\in\bbC\{X\}^*$. Let $N(y)$ be the minimal possible $N$ such that
$y\in\bbC\{X^{1/N}\}$. Suppose that $0<\ord1\,y<\infty$. We write
$$
y(X)=a_1X^{v_1/N}+a_2X^{v_2/N}+\dots\quad a_1,a_2,\dots\neq 0\;,
$$
$0<v_1<v_2<\mbox{}$ integers, $\GCD(N,v_1,v_2,\dots)=1$.
We put $\cycle(y)=\langle y_0,\dots,y_{N-1}\rangle$, where
$$
y_i(X)=a_1\varepsilon^{v_1i}X^{v_1/N}+
       a_2\varepsilon^{v_2i}X^{v_2/N}+\dots,\;
       i=0,\dots,N-1\;,
$$
$\varepsilon$ is a primitive root of 1 of degree $N$.
For $y=0$ we put $\cycle(y)=\langle 0\rangle$. The product
\begin{equation}\label{prod}
[y]:=\prod_{i=0}^{N-1}(Y-y_i(X))\in\bbC\{X,Y\}
\end{equation}
defines a branch. We have $[0]=Y$.
On the other hand, every branch coprime with $X$ can be written in the form
of~(\ref{prod}) up to an invertible factor from $\bbC\{X,Y\}$.

To every $y\in\bbC\{X\}^*$, $\ord1\,y>0$, we assign a {\it generalized
characteristic sequence\/} $(b_0,\dots,b_{\h})$ such that $b_0=N(y)$.
If $N(y)=1$ then $\h(y)=0$ and $(b_0,\dots,b_{\h})=(1)$.
If $N(y)>1$ then we define characteristic positions $(j_1,\dots,j_{\h})$ as
\begin{eqnarray*}
j_k & = & \min\{j>j_{k-1}:\,\GCD(v_0,\dots,v_{j-1})>\GCD(v_0,\dots,v_j)\}\;,
\quad k=1,2,\dots\\
j_{\h} & = & \min\{j:\,\GCD(v_0,\dots,v_j)=1\}.
\end{eqnarray*}
with conventions $j_0=0$ and $v_0=N(y)$. Then
$(b_0,\dots,b_{\h})=(N(y),v_{j_1},\dots,v_{j_{\h}})$.

Let $f=[y]$. We can reconstruct the contact $d(f,X)$ by using $(b_0,\dots,b_{\h})$.
We have $b_0=(f,X)_0$. If $\h(y)=0$ the $d(f,X)=b_0$. If $\h(y)>1$ then
$\ord1\,f=\min\{b_0,b_1\}$. Hence $d(f,X)=b_0/\min\{b_0,b_1\}$. Therefore,
if $b_0<b_1$ then $d(f,X)=1$ ($f,X$ are transverse). If $b_0>b_1$ then
$d(f,X)=b_0/b_1$ ($f,X$ are tangent). If $d(f,X)\notin\bbZ$ then $f,X$
are in the maximal contact.

Let $f_y=f(X,y(X)+Y)\in\bbC\{X^*,Y\}$.
The Newton diagram $\Delta(f_y)$ can be described in terms of characteristics
$(b_0,\dots,b_{\h})$. Let $e_k:=\GCD(b_0,\dots,b_k)$, $k=0,\dots,\h$.
\begin{proper}\label{proper31}
{\rm(see~\cite{GP1991}, Property~3.1 or~\cite{Len2004}, Section~5)}
$$
\Delta(f_{y})=\sum_{k=1}^{\h}
\left\{\Teisss{(b_k/b_0)(e_{k-1}-e_k)}{e_{k-1}-e_k}{20}{10}\right\}
+\left\{\Teisss{\infty}{1}{4}{2}\right\}\;.
$$
\end{proper}
The sum over the empty set equals the zero element $\left\{\Teis{0}{0}\right\}$.
The diagram $\Delta(f_{y_i})$ does not depend on the choice of $y_i\in\cycle(y)$.
For $\Delta\subset\bbR_+^2$, $c>0$ let
$c\Delta=\{ca:\,a\in\Delta\}$.
We define the {\it characteristic Newton diagram\/} of $f$ (with respect
to $X$) as $\Delta^{\chars}_Xf:=(1/\ord1\,f)\Delta(f_y)$.
$$
\epsfbox{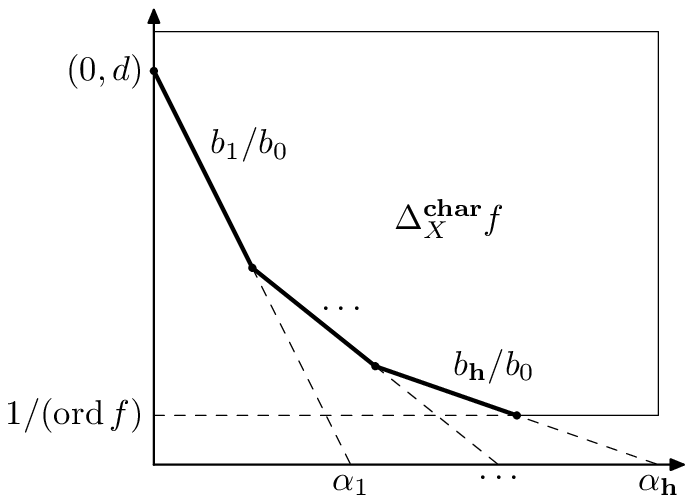}
$$
The diagram $\Delta^{\chars}_Xf$ has a vertex $(0,d)$ on the vertical axis;
$d=d(f,X)$. The distance between the diagram and the
horizontal axis equals $1/(\ord1\,f)$. The inclinations of successive faces
are $b_1/b_0,\dots,b_{\h}/b_0$. Let us denote by
$\alpha_1,\dots,\alpha_{\h}$ the abscissae of points
where the lines determined by succesive faces intersect the horizontal axis.
We restore the characteristic $\char(f)=\{d_1,\dots,d_{\g}\}$
by the following formula, which is a consequence of the Abhyankar
inverse rule~\cite{Campillo}. The characteristic $\{d_1,\dots,d_{\g}\}$ equals
$$
\renewcommand{\arraystretch}{1.2}
\begin{array}{ll}
\{\alpha_1,\dots,\alpha_{\h}\} & \mbox{ if } d(f,X)=1,\;\g=\h\,, \\
\{d\alpha_1,\dots,d\alpha_{\h}\} & \mbox{ if } d(f,X)>1 \mbox{ the contact is maximal, }\g=\h\,,\\
\{d\alpha_2,\dots,d\alpha_{\h}\} & \mbox{ if } d(f,X)>1 \mbox{ the contact is not maximal, }\g=\h-1\,.\\
\end{array}
$$
The sequence $(n_1,\dots,n_{\g})$ can be restored by using
(\ref{n}). We check that
\begin{equation}\label{eq_Ny}
N(y)=d(f,X)n_1\dots n_{\g}\;.
\end{equation}
For a Newton diagram $\Delta$ we define the number $\alpha=\alpha(\varkappa,\Delta)$
which equals the abscissa of the point
where the line of inclination $\varkappa>0$, supporting $\Delta$, intersects
the horizontal axis.
$$
\epsfbox{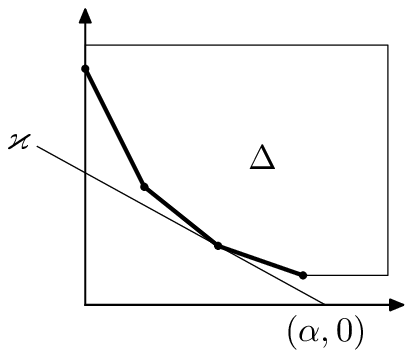}
$$
We also need the inverse operation. We define
the number $\varkappa=\varkappa(\alpha,\Delta)$ as the inclination
of the line supporting $\Delta$ which intersects the horizontal axis
at the point $(\alpha,0)$.

Let us consider $z(X)\in\bbC\{X\}^*$, $\ord1\,z>0$.
Following~\cite{GP1991} let us put
$$
  o_f(z)=\begin{array}[t]{c}\max\\{\scriptstyle y_i\in\cycles(y)}\end{array}
  \ord1(y_i-z)\;.
$$
The number $o_f(z_j)$ does not depend on the choice of $z_j\in\cycle(z)$.
Let $g=[z]$. We have
\begin{equation}\label{eq_contact_exponent}
o_f(z)=o_g(y)=\begin{array}[t]{c}\max\\{\scriptstyle i,j}\end{array}\,
\ord1(y_i-z_j)\;.
\end{equation}
Let us denote (\ref{eq_contact_exponent}) by $\varkappa(f,g,X)$. Let us
consider the ball $B=B(f,g)$, let $\cex(B)$ be the number defined
in (\ref{eq_cex}) and let $d_X(f,g):=\min\{d(f,X),d(g,X)\}$. By using the
inverse rule of Abhyankar we check
\begin{prop}\label{prop_inverse}\mbox{}
\begin{itemize}
\item[\rm(a)] If $d(f,g)=d_X(f,g)$ then
$\displaystyle
 \varkappa(f,g,X)=\min\left\{\frac{1}{d(f,X)},\frac{1}{d(g,X)}\right\}\;,
$
\item[\rm(b)] if $d(f,g)>d_X(f,g)$ then
$\displaystyle
\varkappa(f,g,X)=\frac{\cex(B)-d_X(f,g)+1}{d_X(f,g)}\;.
$
\end{itemize}
\end{prop}
The following property is crucial for our purposes.
\begin{proper}\label{proper33}
{\rm(compare~\cite{GP1991}, Property~3.3)}\\
With previous notation
let $\varkappa=\varkappa(f,g,X)$. Then
$$
d(f,g)=d(X,g)\,\alpha(\varkappa,\Delta^{\chars}_Xf)=
       d(f,X)\,\alpha(\varkappa,\Delta^{\chars}_Xg)\;.
$$
\end{proper}

\vspace{1ex}\noindent{\bf Description of balls in coordinates $X,Y$}\\
Below we characterize an arbitrary ball $B\subset\cB$
in the fixed coordinates $X,Y$. First, we need a fact concerning
Puiseux polynomials. Suppose that $\varphi\in\bbC[X]^*\setminus\{0\}$,
$\ord1\,\varphi>0$. Let
$R_1(\varphi):=d([\varphi],X)\alpha(\deg\varphi,\Delta_X^{\chars}[\varphi])$.
\begin{proper}\label{proper_continuation}
If $f\in\cB$ satisfies $d(f,[\varphi])>R_1(\varphi)$ then there exists
$y\in\Zer\,f$ that is a continuation of $\varphi$.
\end{proper}
\begin{cor}\label{cor_continuation}
With the above assumptions $\char([\varphi])\subset\char(f)$.
\end{cor}
Let $B\subset\cB$ be an arbitrary ball. We measure positions of $B$
with respect to the axes $X$ and $Y$ by the numbers $d(B,X)$
and $d(B,Y)$ (see (\ref{eq_dBh})). In any case we have $d(B,X)\leq d(B)$
and $d(B,Y)\leq d(B)$.
\begin{proper}\label{proper_balls1}
{\rm(classification of balls)}\\
Let $R=d(B)$. Then:
\begin{itemize}
\item[\rm(I)] If $R=1$ then $B=\cB$.
\item[\rm(II)] If $1<R=d(B,X)$ then $B=B(X,R)$\quad
\item[\rm(III)] If $1<R=d(B,Y)$ then $B=B(Y,R)$.
\item[\rm(IV)] If $R>\max\{d(B,X),d(B,Y)\}$ then there exists
$\varphi\in\bbC[X]^*\setminus\{0\}$, $\ord1\,\varphi>0$ such that
$R_1(\varphi)<R$ and $B=B([\varphi],R)$.
\end{itemize}
Moreover, the classes of balls described above are disjoint.
\end{proper}
A ball $B$ belongs to the joined class (III)+(IV) if and only if
$d(B,X)<d(B)$. Since for $\varphi=0$ we have $[\varphi]=Y$ it is convenient
to put $R_1(0)=1$. After that in the joined class (III)+(IV) we have
the following
\begin{proper}\label{rem_uniqueness}{\rm
Let $B$ be a ball satifying $d(B,X)<d(B)$. Then if
$B=B([\varphi],R)=B([\tilde\varphi],R)$ where
$\varphi,\tilde\varphi\in\bbC[X]^*$, $\ord1\,\varphi>0$,
$\ord1\,\tilde\varphi>0$, $R_1(\varphi)<R$, $R_1(\tilde\varphi)<R$ then
$[\varphi]=[\tilde\varphi]$ (i.e. $\varphi$ and $\tilde\varphi$ are in
the same cycle).
}
\end{proper}
First, we compute the invariants $\char(B)$, $\nu(B)$, $n(B)$ which depend
only on the ball.
Let $B$ be a ball from Property~\ref{proper_balls1} with the radius
$R=d(B)$ and let $d_X=d(B,X)$. We assign to $B$ the numbers
\begin{equation}\label{eq_kappa}
\renewcommand{\arraystretch}{1.2}
\varkappa=\left\{\begin{array}{ll}
1 & \mbox{ in case (I)}\\
1/R & \mbox{ in case (II)}\\
R & \mbox{ in case (III)}\\
\varkappa(R/d_X,\Delta_X^{\chars}[\varphi]) & \mbox{ in case (IV)}\\
\end{array}\right.
\end{equation}
and
\begin{equation}\label{eq_N}
\renewcommand{\arraystretch}{1.2}
N=\left\{\begin{array}{ll}
1 & \mbox{ in cases (I), (II), (III)}\\
N(\varphi) & \mbox{ in case (IV)}.\\
\end{array}\right.
\end{equation}
Let us write $\varkappa=m/(N\nn)$, $\GCD(\nn,m)=1$.
\begin{prop}\label{prop_balls2}
Let $B$ be a ball from Property~\ref{proper_balls1}. Then with the previous
notation we have
\begin{center}
\renewcommand{\arraystretch}{1.2}
\begin{tabular}{l|c|c|c}
{\rm case} & $\char(B)$ & $\nu(B)$ & $n(B)$ \\
\hline
{\rm(I)} & $\varnothing$ & $1$ & $1$ \\
{\rm(II)} & $\varnothing$ & $1$ & $m$ \\
{\rm(III)} & $\varnothing$ & $1$ & $\nn$ \\
{\rm(IV)} & $\char([\varphi])$ & $N/d_X$ & $\nn$ \\
\end{tabular}
\end{center}
\end{prop}
Proof. (I), (II), (III) follow directly from the definition.

\vspace{1ex}\noindent(IV).
Since $[\varphi]$ is a center of $B$ we have
$\char([\varphi])\supset\char(B)$. Let $f\in\cB$. We have
$d(f,[\varphi])\geq R>R_1(\varphi)$. Then we use
Corollary~\ref{cor_continuation}. We obtain $\nu(B)=N/d_X$ from (\ref{eq_Ny}).
To show $n(B)=\nn$ we use (\ref{nb}) and (\ref{eq_cex}). We check that
$$
\left\{n\geq 1:\,d(B)\in\frac{\bbN}{\nu(B)^2n}\right\}=
\left\{n\geq 1:\,\cex(B)\in\frac{\bbN}{\nu(B)n}\right\}\;.
$$
We finish by using Proposition~\ref{prop_inverse}~\rule{1ex}{1ex}

Now, for a ball from family (\ref{def_T}) determined by the germ
we want to compute the numbers $t^{(1)}$ and $t^{(2)}$ of direct successors
in the Eggers tree. These numbers depend not only on the ball but also on the
germ. For $f\in\bbC\{X^{1/N},Y\}$ by $r^{(N)}(f)$ we denote the number of
pairwise coprime factors of $f$ in $\bbC\{X^{1/N},Y\}$;
$r_0^{(N)}(f)$ stands for the number of factors different from $X$ and $Y$.
For $N=1$ we write $r(f)$ and $r_0(f)$, respectively. We put
$\varepsilon_X(f)=1$ if $X$ appears as a factor of $f$ and
$\varepsilon_X(f)=0$, otherwise. Analogously we define $\varepsilon_Y(f)$.
Let $(a,b)$ be a vector ($a,b>0$) and let
$f=\sum c_{\alpha\beta}X^\alpha Y^\beta$. We define the initial form
$\In_{(a,b)}f$ of $f$ with respect to $(a,b)$ as
$\sum c_{\alpha\beta}X^\alpha Y^\beta$ where $(\alpha,\beta)\in\supp\,f$ and
$\alpha\,a+\beta\,b=\inf\{\alpha\,a+\beta\,b:\,(\alpha,\beta)\in\supp\,f\}$.
For $S\in\cN_f$ we have $\In(f,S)=\In_{(|S|_2,|S|_1)}f$. For a generic
$(a,b)$ $\In_{(a,b)}f$ is a monomial.

For a pair $(\varphi,\varkappa)$, $\varphi\in\bbC[X]^*$, $\ord1\,\varphi>0$,
$\varkappa>\deg\varphi$ we consider two Puiseux series
$y=\varphi+aX^\varkappa+\dots$,
$z=\varphi+bX^\varkappa+\dots$, where $a,b\in\bbC$. Let $N=N(\varphi)$,
$\varkappa=m/(N\nn)$, $\GCD(\nn,m)=1$. Let $f=[y]$ and $g=[z]$.
We need a tool for estimating the contact $d(f,g)$. Let
\begin{equation}\label{eq_R}
\renewcommand{\arraystretch}{1.2}
R_{\varphi,\varkappa}=\left\{\begin{array}{ll}
\max\{\varkappa,1/\varkappa\} & \mbox{if }\varphi=0\;,\\
d([\varphi],X)\,\alpha(\varkappa,\Delta_X^{\chars}[\varphi])
& \mbox{if }\varphi\neq 0\;.
\end{array}\right.
\end{equation}
\begin{proper}\label{proper_2series}\mbox{}
Assume that if $\varphi=0$ and $\varkappa<1$ then $a,b$ are nonzero.
Otherwise, $a,b$ are arbitrary. Then $d(f,g)\geq R_{\varphi,\varkappa}$ and
the following conditions are equivalent:
\begin{itemize}
\item[\rm(1)] $d(f,g)>R_{\varphi,\varkappa}$,
\item[\rm(2)] $a^{\nn}=b^{\nn}$,
\item[\rm(3)] $aX^\varkappa$ is a root of $\In_{(1,\varkappa)}g_\varphi$ in $\bbC\{X\}^*$,
\item[\rm(4)] $bX^\varkappa$ is a root of $\In_{(1,\varkappa)}f_\varphi$ in $\bbC\{X\}^*$,
\item[\rm(5)] $\In_{(1,\varkappa)}f_\varphi$ and $\In_{(1,\varkappa)}g_\varphi$
can be written us $X^{\zeta_i}(Y^{\nn}-cX^{m/N})^{\eta_i}$ up to nonzero constants
with $\zeta_i\geq 0,\eta_i>0$ $(i=1,2)$ and $c=a^{\nn}=b^{\nn}$ {\rm(}possibly zero{\rm)}.
\end{itemize}
\end{proper}
Now, let us consider a germ $f=0$, $f\in\bbC\{X,Y\}$ reduced,
and an arbitrary ball $B$ from Property~\ref{proper_balls1}. The formulas for
computing the numbers $t^{(1)}_f(B)$ and $t^{(2)}_f(B)$ are presented in the
following proposition which is a consequence of Property~\ref{proper_2series}.
If $B$ is characteristic ($n(B)>1$) then in order to determine $t^{(1)}_f(B)$
we consider the equivalency class of a branch $h\in B$ such that
$d(B)\notin\char(h)$ (see proof of Proposition~\ref{prop_t}(c)). In cases
(II), (III), (IV) as $h$ we choose $X$, $Y$ and $[\varphi]$, respectively.
\begin{prop}\label{prop_table}\mbox{}{\rm
\begin{center}
\renewcommand{\arraystretch}{1.2}
\begin{tabular}{ll|l|l}
case & & $t_f^{(1)}(B)$ & $t_f^{(2)}(B)$ \\
\hline
(I)  & & $r(\In\,f)$  & 0 \\
(II) & $n(B)=1$ & $r_0(\In_{(R,1)}f)+\varepsilon_X(\In_{(R,1)}f)$ & $0$\\
     & $n(B)>1$ & $\varepsilon_X(\In_{(R,1)}f)$ & $r_0(\In_{(R,1)}f)$ \\
(III) & $n(B)=1$ & $r_0(\In_{(1,R)}f)+\varepsilon_Y(\In_{(1,R)}f)$ & $0$\\
     & $n(B)>1$ & $\varepsilon_Y(\In_{(1,R)}f)$ & $r_0(\In_{(1,R)}f)$ \\
(IV) & $n(B)=1$ & $r_0^{(N)}(\In_{(1,\varkappa)}f_\varphi)+\varepsilon_Y(\In_{(1,\varkappa)}f_\varphi)$ & $0$\\
     & $n(B)>1$ & $\varepsilon_Y(\In_{(1,\varkappa)}f_\varphi)$ & $r_0^{(N)}(\In_{(1,\varkappa)}f_\varphi)$ \\
\end{tabular}
\end{center}
} %
\end{prop}

\vspace{1ex}\noindent{\bf Proof of Theorem~\ref{thm_factorization}}\\
We consider a singularity $f=0$, $f\in\bbC\{X,Y\}$ reduced, and
a regular parameter $\lambda$. Without loss of generality we can
assume that $\lambda=X$. Then $\bJ(\lambda,f)=\P{f}/\P{Y}$.
For every $\varphi\in\Track(f)$ we consider the Newton polygon
$\cN_\varphi:=\cN^{\deg\varphi}(f_\varphi)$. For $\varphi=0$
we obtain the classical Newton polygon $\cN_f$. Let $\cB_f=\{f_1,\dots,f_r\}$
be the set of branches of $f$. We can write
$f=X^{\delta_X(f)}\tilde{f}$ where $(\tilde{f},X)_0=p>0$. Clearly
$\Zer(f)=\Zer(\tilde{f})$ and $\Zer(\P{f}/\P{Y})=\Zer(\P{\tilde f}/\P{Y})$.
We apply Theorem~\ref{thm_Len2004} to $\tilde{f}$.
As in (\ref{factorizations}) we consider the factorization
$\P{f}/\P{Y}=X^{\delta_X(f)}g_1\dots g_u$ where the branches
$g_1,\dots,g_u$ are coprime with $X$. We have
$\Zer(\P{f}/\P{Y})=\Zer\,g_1\oplus\dots\oplus\Zer\,g_u$.
With the notation of Theorem~\ref{thm_Len2004} we state the following
\begin{lm}\label{lm_balls} Let $z(X)\in\Zer(\P{f}/\P{Y})$. Then
\begin{itemize}
\item[\rm(I)] if $z\in\bA_S$, $S\in\cN_f$, $\incl(S)=1$ then $B_f([z])=\cB$;
\item[\rm(II)] if $z\in\bA_S$, $S\in\cN_f$, $\incl(S)<1$ then
\begin{itemize}
\item[$(\cdot)$]
if $z$ is of the first kind then $B_f([z])=B(X,1/\incl(S))$;
\item[\rm$(\cdot\cdot)$]
if $z$ is of the second kind then $B_f([z])<B(X,1/\incl(S))$;
\end{itemize}
\item[\rm(III)] if $z\in\bA_S$, $S\in\cN_f$, $\incl(S)>1$ then $B_f([z])=B(Y,\incl(S))$;
\item[\rm(IV)] if $z\in\bA_{\varphi,S}$, $\varphi\neq 0$, $S\in\cN_\varphi$
then $B_f([z])=B([\varphi],R_{\varphi,S})$ where
$R_{\varphi,S}:=R_{\varphi,\incls(S)}$.
\end{itemize}
\end{lm}
Proof. By Corollary~\ref{cor_useful_balls}(b) we have
\begin{equation}\label{eq_ball}
\renewcommand{\arraystretch}{1.2}
B_f([z])=B(h,R)\Leftrightarrow
\left(\begin{array}{c}\max\{d(f_1,[z]),\dots,d(f_r,[z])\}=R\\
\mbox{ and }d([z],h)\geq R\end{array}\right)\;.
\end{equation}
Recall that
\begin{equation}\label{eq_ord}
\ord1\,z=\frac{d([z],Y)}{d([z],X)}\;.
\end{equation}

\vspace{1ex}\noindent(I). Let $z\in\bA_S$, $\incl(S)=1$. Let
$I_S=\{i\in\{1,\dots,r\}:\,d(f_i,X)=d(f_i,Y)=1\}$. Assume that $z$ is of the
first kind (0-first kind). That is $\ord1\,z=\incl(S)=1$ (the styles
of edges are not expressed in the sketches).
$$
\epsfbox{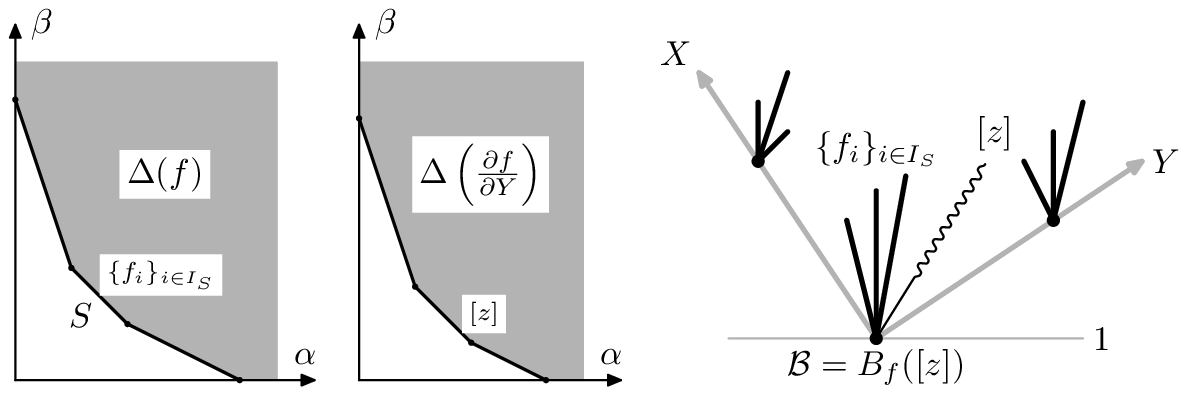}
$$
Hence $d([z],X)=d([z],Y)=1$ and by the definition of $\bA_S$ $\In\,z$ is not a
root of $\In(f,S)=\In f=\prod_{i=1}^r\In\,f_i$. For $i\in I_S$ there exists
$S_i\in\cN(f_i)$ parallel to $S$ and $\In\,f_i=\In(f_i,S_i)$. Therefore,
$d(f_i,[z])=1$ by Property~\ref{proper_2series}. For $i\notin I_S$ we
have $d(f_i,X)>1$ or $d(f_i,Y)>1$. Hence $d(f_i,[z])=1$ by ($D_3'$) and we
have $B_f([z])=\cB$ by (\ref{eq_ball}).

If $z$ is of the second kind then $S$ touches the horizontal axis
(Proposition~\ref{prop_LenMast2009}(ii)).
$$
\epsfbox{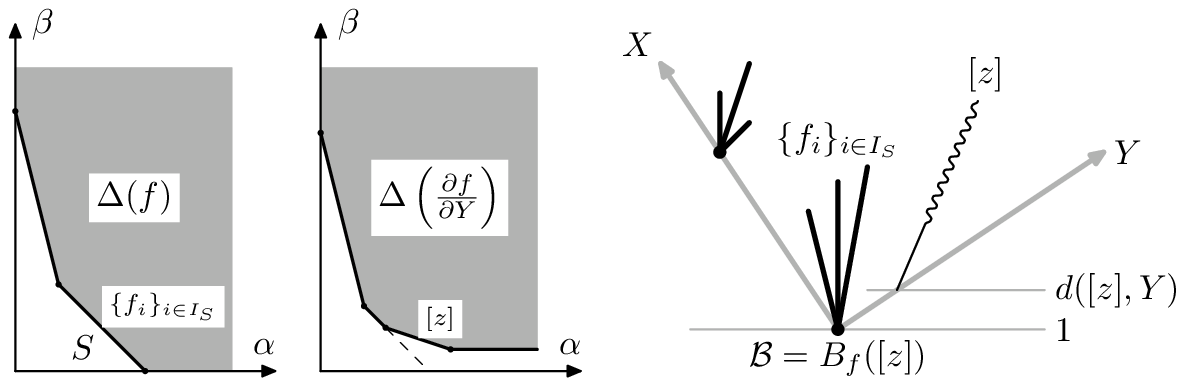}
$$
We have $d(f_i,Y)=1$ for $i\in\{1,\dots,r\}$ and $\ord1\,z>\incl(S)=1$ by
Proposition~\ref{prop_LenMast2009}(iii). Hence $d([z],Y)>1$ by (\ref{eq_ord}).
From ($D_3'$) we obtain $d(f_i,[z])=1$ for $i\in\{1,\dots,r\}$ which gives
$B_f([z])=\cB$ by~(\ref{eq_ball}).

\vspace{1ex}\noindent (II). Let $z\in\bA_S$, $\incl(S)<1$.
Let $I_S=\{i\in\{1,\dots,r\}:\,d(f_i,X)=1/\incl(S)\}$.
Assume that $z$ is of the first kind.
$$
\epsfbox{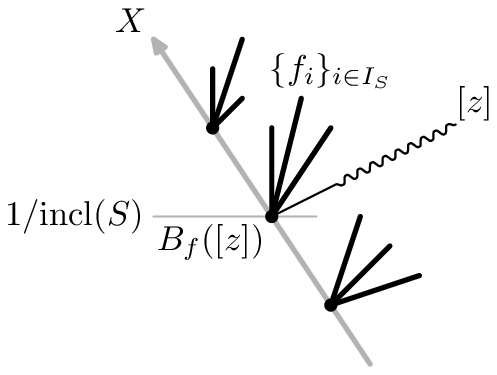}
$$
Hence $d([z],X)=1/\incl(S)$ and $\In z$ is not a root of
$\In(f,S)=\prod_{i=1}^r\In_{\vec v}f_i$ with $\vec v=(1,\incl(S))$.
Analogously, as earlier by using Property~\ref{proper_2series} and ($D_3'$)
we conclude that $B_f([z])=B(X,1/\incl(S))$.

If $z$ is of the second kind, then $S$ touches the horizontal
axis and $\ord1\,z>\incl(S)$ by Proposition~\ref{prop_LenMast2009}(iii).
$$
\epsfbox{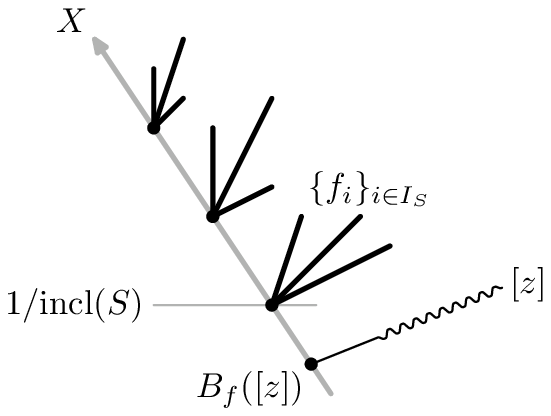}
$$
By using (\ref{eq_ord}) we obtain $d([z],X)<1/\incl(S)$. For
$i\in\{1,\dots,r\}$ the inequality $d(f_i,X)\geq 1/\incl(S)$ results in
$d(f_i,[z])=d([z],X)$. Hence $B_f([z])=B(X,d([z],X))<B(X,1/\incl(S))$.

\vspace{1ex}\noindent(III)+(IV).
Let $z\in\bA_{\varphi,S}$, $S\in\cN_\varphi$ (if $\varphi=0$ then $\incl(S)>1$).
Let $I_{\varphi,S}=\{i\in\{1,\dots,r\}:\,d(f_i,[\varphi])=R_{\varphi,S}\}$.
Assume that $z$ is of the $\varphi$-first kind.
$$
\epsfbox{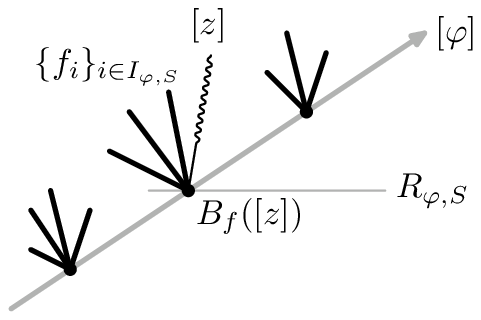}
$$
Then $\ord1(z-\varphi)=\incl(S)$ and $d([z],[\varphi])=R_{\varphi,S}$
by Property~\ref{proper_2series}. By the definition of $\bA_{\varphi,S}$
$\In(z-\varphi)$ is not a root of
$\In(f_\varphi,S)=\prod_{i=1}^r\In_{\vec{v}}(f_i)_\varphi$
with $\vec{v}=(1,\incl(S))$. As in (I) by using
Property~\ref{proper_2series} and ($D_3'$) we conclude that
$B_f([z])=B([\varphi],R_{\varphi,S})$.

If $z$ is of the $\varphi$-second kind, then $S$ touches the horizontal axis,
$\ord1(z-\varphi)>\incl(S)$ and $d(f_i,[\varphi])\leq R_{\varphi,S}$
for $i\in\{1,\dots,r\}$ by Proposition~\ref{prop_LenMast2009}.
$$
\epsfbox{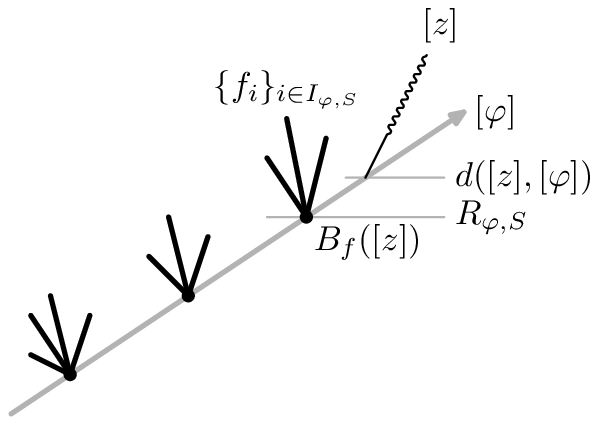}
$$
By Proposition~\ref{proper_2series} $d([z],[\varphi])>R_{\varphi,S}$.
Therefore, by ($D_3'$) $d(f_i,[z])=d(f_i,[\varphi])$.
Hence $d(f_i,[z])\leq R_{\varphi,S}$ with equality for $i\in I_{\varphi,S}$.
By using (\ref{eq_ball}) we obtain
$B_f([z])=B([\varphi],R_{\varphi,S})$~\rule{1ex}{1ex}

\vspace{1ex}\noindent
Let $\Track(f)^*=\Track(f)\setminus\{0\}$. Following Lemma~\ref{lm_balls}
let us define three sets of balls. We denote by $\cT_1$ the set of
balls $B(X,1/\incl(S))$ where $S\in\cN_f$ with $\incl(S)<1$, by $\cT_2$ the
set of balls $B(Y,\incl(S))$ where $S\in\cN_f$ with $\incl(S)\geq 1$
and by $\cT_3$ the set of balls $B([\varphi],R_{\varphi,S})$ where
$\varphi\in\Track(f)^*$, $S\in\cN_\varphi$.

The chain $\oK_f(X)=\{B(f_1,X),\dots,B(f_r,X)\}$ plays an important role in our
approach (see Property~\ref{proper_extension} and the text before).
\begin{proper}\label{proper_chainK}
\renewcommand{\arraystretch}{1.2}
{\rm(description of the chain $\oK_f(X)$)}
\begin{itemize}
\item[\rm(a)] $\cT_1\subset\oK_f(X)\subset\{\cB\}\cup\cT_1\cup\{B(X,X)\}$ with
\begin{itemize}
\item $\cB\in\oK_f(X)\Leftrightarrow\max\,\incl(\cN_f)\geq 1$,
\item $B(X,X)\in\oK_f(X)\Leftrightarrow\delta_X(f)>0$.
\end{itemize}

\item[\rm(b)] $\displaystyle
\min\,\oK_f(X)=\left\{\begin{array}{lcl}B(X,1/\max\,\incl(\cN_f)) &
\mbox{ if } & \max\,\incl(\cN_f)<1\\
\cB & \mbox{ if } & \max\,\incl(\cN_f)\geq 1\\
\end{array}\right.$

\item[\rm(c)] $B_f(X)=\left\{\begin{array}{lcl}
B(X,X) & \mbox{ if } & \delta_X(f)>0\\
B(X,1/\min\,\incl(\cN_f)) & \mbox{ if } & \delta_X(f)=0\mbox{ and }%
\min\,\incl(\cN_f)<1\\
\cB & \mbox{ if } & \delta_X(f)=0\mbox{ and }\min\,\incl(\cN_f)\geq 1\\
\end{array}\right.$
\end{itemize}
\end{proper}
Lemma~\ref{lm_balls} and Property~\ref{proper_chainK} allow us to finish
the proof of part (ii) of the factorization theorem. We want to show
the equality $d(g,X)=d_f(g,X)$ for every branch $g$ of $\P{f}/\P{Y}$ where
$d_f(g,X)$ is defined by~(\ref{contact_rel}). If $B_f(g)\neq B_f(X)$ then
the equality is a direct consequence of ($D_1$--$D_3$). Assume that
$B_f(g)=B_f(X)$. To finish the proof it suffices to find a branch
$f_{i_0}\in\cB_f$ such that
\begin{equation}\label{eq_fi0}
d(g,X)=\min\{d(f_{i_0},g),d(f_{i_0},X)\}\;.
\end{equation}
Let us consider the cases determined for $B_f(X)$ by
Property~\ref{proper_chainK}.

\vspace{1ex}\noindent If $\delta_X(f)>0$ then $B_f(X)=B(X,X)$.
The equality $B_f(g)=B_f(X)$ leads to $g=X$. With $f_{i_0}=X$
we obtain $\infty$ on both sides of (\ref{eq_fi0}).

\vspace{1ex}\noindent If $\delta_X(f)=0$ and $\min\,\incl(\cN_f)<1$ then
$B_f(X)=B(X,1/\incl(S))$ where $S$ has the minimal inclination.
Now, every branch of $\P{f}/\P{Y}$ has the form $g=[z]$ where
$z\in\Zer(\P{f}/\P{Y})$. From $B_f([z])=B(X,1/\incl(S))$ and
Lemma~\ref{lm_balls} it follows that $z\in\bA_S$ and $z$ is of the
first kind. We obtain (\ref{eq_fi0}) for every $i_0\in I_S$ where
$I_S$ is defined in the proof of Lemma~\ref{lm_balls}.

\vspace{1ex}\noindent If $\delta_X(f)=0$ and $\min\,\incl(\cN_f)\geq 1$ then
$B_f(X)=\cB$. Assume that $B_f([z])=B_f(X)=\cB$. The case of
Lemma~\ref{lm_balls}(II)($\cdot\cdot$) is impossible. The only
possibility is that $z\in\bA_S$, $S\in\cN_f$, $\incl(S)=1$. As in the
proof of Lemma~\ref{lm_balls} we show that (\ref{eq_fi0}) is satisfied
for every $i_0\in\{1,\dots,r\}$~\rule{1ex}{1ex}

We are in a good position to finish the proof of parts (i) and (iii) of
Theorem~\ref{thm_factorization}.
Let $\cT=\cT_1\cup\cT_2\cup\cT_3$ and let $\cT':=\cT\cup\{\cB\}$.
By Property~\ref{proper_chainK} $\min\,\oK_f(X)\in\cT'$.
Let $\Zer(\P{f}/\P{Y})=\{z_1,\dots,z_{p-1}\}$ as earlier. Let $B\in\cT'$.
For $B\neq\min\,\oK_f(X)$ we put
$J_B=\{j\in\{1,\dots,p-1\}:\,B_f([z_j])=B\}$ and for $B=\min\,\oK_f(X)$
we put $J_B=\{j\in\{1,\dots,p-1\}:\,B_f([z_j])\leq B\}$.
For every $B\in\cT'$ let us define $h_B=\prod_{j\in J_B}(Y-z_j)$ with
convention $\prod_\varnothing=1$. Clearly $h_B\in\bbC\{X,Y\}$.
By Lemma~\ref{lm_balls} and Theorem~\ref{thm_Len2004} for $B\in\cT$
the number $(h_B,X)_0$ equals
$$
\renewcommand{\arraystretch}{1.2}
\begin{array}{lcl}
t(f,S)-1 & \mbox{ if } & B=B(X,1/\incl(S)),\,S\in\cN_f,\,\incl(S)<1\\
t(f,S)-1 & \mbox{ if } & B=B(Y,\incl(S)),\,S\in\cN_f,\,\incl(S)\geq 1\\
N(\varphi)(t(f_\varphi,S)-1) & \mbox{ if } &
\varphi\in\Track(f)^*,\,S\in\cN_\varphi\\
\end{array}
$$
For every $B\in\cT'$ we put
$$
\renewcommand{\arraystretch}{1.2}
m'(B)=\left\{\begin{array}{ccl}
(h_B,X)_0 & \mbox{ if } & B\in\cT\\
0 & \mbox{ if } & B\notin\cT\;.\\
\end{array}\right.
$$
The second case is possible only if $B=\cB$ and $\cB\notin\cT$.

To finish the proof of Theorem~\ref{thm_factorization}(iii)
we need the following
\begin{prop}\label{prop_multiplicities}
For $B\in\cT'$
the number $m'(B)$ equals:\\
$d(B,X)\nu(B)[t_f^{(1)}(B)+n(B)t_f^{(2)}(B)-1]\rule{0ex}{3ex}$
for $B\in\cT'\setminus\oK_f(X)$;\\
$d(B)n(B)[t_f(B)-1+\sigma^{\maxs}_{f,X}(B)]-\sigma^{\mins}_{f,X}(B)\rule{0ex}{3ex}$
for $B\in\oK_f(X)$.
\end{prop}
Proof. Assume first that $B\in\cT_3$ (hence $B\notin\oK_f(X)$ by
Property~\ref{proper_balls1}). There
exist $\varphi\in\Track(f)^*$ and $S\in\cN_\varphi$ such that
$B=B([\varphi],R_{\varphi,S})$. We have $\incl(S)=m/(N\nn)$, $\GCD(\nn,m)=1$,
$N=N(\varphi)$, $\nn=n(B)$ (Proposition~\ref{prop_balls2}). Let us observe that
$t(f_\varphi,S)=n(B)r_0^{(N)}(\In(f_\varphi,S))+%
\varepsilon_Y(\In(f_\varphi,S))$. We finish by using (\ref{eq_Ny}),
Property~\ref{dBh} and Property~\ref{prop_table}.

\vspace{1ex}\noindent If $B\in\cT_2\setminus\{\cB\}$ (hence
$B\notin\oK_f(X)$) we put $\varphi=0$ above.

\vspace{1ex}\noindent If $B\in\cT_1$ then $B\in\oK_f(X)$. We have
$B=B(X,1/\incl(S))$, $S\in\cN_f$, $\incl(S)<1$. In this case
$\sigma^{\maxs}_{f,X}(B)=1-\varepsilon_X(\In(f,S))$ and
$\sigma^{\mins}_{f,X}(B)=1-\varepsilon_Y(\In(f,S))$. We have
$t(f,S)=d(B)n(B)r_0(\In(f,S))+\varepsilon_Y(\In(f,S))$. We finish by
using Proposition~\ref{prop_table}.

\vspace{1ex}\noindent Now, consider the case when $\cB\in\cT$. Hence there
exists $S\in\cN_f$ with $\incl(S)=1$ and $\cB\in\oK_f(X)$. In this case
$\sigma^{\maxs}_{f,X}(B)=1-\varepsilon_X(\In(f,S))$ and
$\sigma^{\mins}_{f,X}(B)=1$. We finish by using Proposition~\ref{prop_table}.

\vspace{1ex}\noindent If $\cB\notin\cT$ then there does not exist a face
of inclination~1 in $\cN_f$. We check the appropriate formulas
directly~\rule{1ex}{1ex}

\vspace{2ex}To end the proof of Theorem~\ref{thm_factorization}~(i) and~(iii)
we observe that for $B\in\cT'$ if $m'(B)>0$ then
$B\in\cE(f)\cup\{B_f(X)\}$~\rule{1ex}{1ex}

\section*{Acknowledgements}
The author wishes to express his thanks to Piotr Tworzewski
for his support (from KBN grant)
throughout preliminary versions of this paper
and to Arkadiusz P\l{}oski and to Janusz Gwo\'{z}dziewicz
for their help with preparing its final version.

\vspace{2ex}\noindent
{\sc Department of Mathematics, Kielce University of Technology,\\
AL. 1000 L PP 7, 25-314 Kielce, Poland\\
e-mail: }\verb+ztpal@tu.kielce.pl+

\end{document}